\documentclass[11pt]{article}%
\usepackage{amsmath}
\usepackage{amsfonts}
\usepackage{amssymb}
\usepackage{graphicx}
\usepackage{subfigure}
\usepackage{citesort}

\newtheorem{theorem}{Theorem}

\newtheorem{proposition}[theorem]{Proposition}

\begin{document}

\title
{2D and 3D reconstructions in acousto-electric tomography}

\author{Peter Kuchment and Leonid Kunyansky}

\maketitle

\begin{abstract}
We propose and test stable algorithms for the reconstruction of the internal
conductivity of a biological object using acousto-electric measurements.
Namely, the conventional impedance tomography scheme is supplemented by
scanning the object with acoustic waves that slightly perturb the conductivity
and cause the change in the electric potential measured on the boundary of the
object. These perturbations of the potential are then used as the data for the
reconstruction of the conductivity. The present method does not rely on
``perfectly focused'' acoustic beams. Instead, more realistic propagating
spherical fronts are utilized, and then the measurements that would correspond
to perfect focusing are synthesized. In other words, we use \emph{synthetic
focusing}. Numerical experiments with simulated data show that our techniques
produce high quality images, both in $2D$ and $3D$, and that they remain accurate
in the presence of high-level noise in the data. Local uniqueness and stability
for the problem also hold. \end{abstract}

%

\section*{Introduction}

Electrical Impedance Tomography (EIT) is a harmless and inexpensive imaging
modality, with important clinical and industrial applications. It aims to
reconstruct the internal conductivity of a body using boundary electric
measurements (see, e.g., \cite{Bor02,BB1,CIN,Cipra}). It is well known that,
regretfully, it suffers from inherent low resolution and instability. To
bypass this difficulty, various versions of a new hybrid technique, sometimes
called Acousto-Electric Tomography (AET), have been introduced recently
\cite{Ammari_EIT,KuKuAET,Wang_AET,Cap}. (See also \cite{Scherzer} for a
different way to recover the conductivity using combination of ultrasound and
EIT). AET utilizes the electro-acoustic effect, i.e. occurrence of small
changes in tissue conductivity as the result of applied acoustic pressure
\cite{AE2,AE3}. Although the effect is small, it was shown in \cite{Wang_AET}
that it provides a signal that can be used for imaging the conductivity. It
has been understood \cite{Ammari_EIT,KuKuAET,Cap} that if one could apply
concentrated pressure at a given point inside the body and then measure the
resulting change in impedance measurements, the knowledge of the perturbation
point would have a stabilizing effect on the reconstruction in otherwise
highly unstable EIT. It has been proposed to use a tightly focused ultrasound
beam as a source of such point-like acoustic pressure \cite{Ammari_EIT}.
However, since perfect focusing of acoustic waves is hard to achieve in
practice (see, e.g., \cite{Localized}), an alternative \emph{synthetic
focusing} approach was developed in \cite{KuKuAET}. Namely, the medium is
perturbed by a series of more realistic propagating spherical acoustic fronts
with centers lying outside of the object (other options, e.g. plane waves or
monochromatic spherical waves could also be used \cite{KuKuAET}). The
resulting changes in the values of electric potential on the boundary of the
object are recorded. Then the data that would have been collected, if perfect
focusing were possible, are synthesized mathematically. Such synthesis happens
to be equivalent to the well established inversion in the so called
\emph{thermoacoustic tomography} (see, e.g., the surveys
\cite{Wang_book,MXW_review,KuKu}). Of course, for accurate synthesis the
acoustic properties of the medium should be known. In breast imaging, for example,
the speed of sound in the tissue can be well approximated by a constant,
and application of AET in this area looks very promising. In the inhomogeneous
medium synthetic focusing is possible if its acoustic parameters
are reconstructed beforehand (for example, using methods of ultrasound
tomography). The results of first numerical experiments
presented in \cite{KuKuAET} confirm the feasibility of the synthetic
focusing.

In this article, we describe a stable and efficient local algorithm for the
AET problem. From the formulas we present one can easily infer the local
uniqueness and stability of the reconstruction. However, after this work was
done, the authors have learned of the paper \cite{Cap}, some results of which
(Propositions 2.1, 2.2) imply uniqueness and Lipschitz stability in the
similar setting (see also \cite{bon} for the presentation of such a local
result). We thus address these issues only briefly here.

The presented algorithm involves two steps. First, it synthesizes the data
corresponding to perfectly focused ultrasound perturbations from the data
obtained using more realistic spherical waves. Here the known smallness of the
acousto-electric effect \cite{AE2,AE3,Wang_AET} is crucial, since it permits
linearization with respect to the acoustic perturbation and thus makes
synthetic focusing possible. Second, the algorithm reconstructs the
conductivity from the data corresponding to perfectly focused perturbations.
This second step, from measured data to the conductivity, is non-linear. We
develop a linearized algorithm, assuming that the conductivity is close to a
known one. The numerical examples that we provide show that this approach
works surprisingly well even when the initial guess is very distinct from the
correct conductivity. One can apply iterations for further improvements.

To the best of the authors' knowledge, the first step of our method (synthetic
focusing) has not been discussed previously in works on AET, except for a
brief description in our papers \cite{KuKuAET,Oberw}. On the other hand, three
different approaches to reconstruction using perfectly focused beam (the
second step of our algorithm) have been recently proposed
\cite{Ammari_EIT,KuKuAET,Cap,Oberw}. Let us thus indicate the differences with
these recent works.

In \cite{Ammari_EIT}, two boundary current profiles were used and the
problem of reconstructing the conductivity was reduced to a numerical solution
of a (non-linear) PDE involving the $0$-Laplacian. In \cite{KuKuAET,Oberw}, by
a rather crude approximation, we reduced the reconstruction problem to solving
a transport equation (a single current was used). Unfortunately, in the case
of noisy measurements the errors tend to propagate along characteristics,
producing unpleasant artifacts in the images, which can be reduced by
iterations. There is also a version of this procedure that involves an elliptic
equation and thus works
better. In \cite{Cap}, two current profiles are used in $2D$
(three profiles in $3D$), the problem is reduced to
a minimization problem, which is then solved numerically. In the present paper
we also use two currents in $2D$ (two or three in $3D$) and, on the second step,
we utilize the same data as in \cite{Cap}. Unlike \cite{Cap}, in our work the
reconstruction problem is solved, under the assumption that the conductivity is
close to some initial guess, by a simple algorithm, which even on the first
step produces good images, improved further by iterations. The algorithm
essentially boils down to solving a Poisson equation. Numerical experiments
show high quality reconstructions, quite accurate even in the presence of very
significant noise. Reconstructions remain accurate when the true conductivity
differs significantly from the initial guess.

The rest of the paper is organized as follows: Section \ref{S:formulation}
contains the formulation of the problem. It also addresses the focusing issue.
The next Section \ref{S:reconstruction} describes the reconstruction algorithm,
stability of which is discussed in Section \ref{S:stability}. Numerical
implementation and results of reconstruction from simulated data in $2D$ are
described in Section \ref{S:numerical}. Sections \ref{S:3drec} and \ref{S:3d}
are addressing the $3D$ case. Section \ref{S:remarks} is devoted to final
remarks and conclusions.

\section{Formulation of the problem}\label{S:formulation}

Let $\sigma(x)$ be the conductivity of the medium within a bounded
region~$\Omega$. Then the propagation of the electrical currents through
$\Omega$ is governed by the divergence equation
\begin{equation}
\nabla\cdot\sigma(x)\nabla u(x)=0,x\in\Omega. \label{original}
\end{equation}
or, equivalently%
\begin{equation}
\Delta u(x)+\nabla u(x)\cdot \nabla \ln \sigma (x)=0,  \label{original1}
\end{equation}%
where $u(x)$ is the electric potential. Let us assume that $\sigma -1$ is
compactly supported within region $\Omega ,$ and that $\sigma (x)=1$ in the
neighborhood of the boundary $\partial \Omega .$ We also assume that the
currents $J=\sigma \frac{\partial }{\partial n}u(x)$ through the boundary
are fixed and the values of potential $u$ are measured on the boundary $%
\partial \Omega $.

The acoustic wave propagating through the object slightly perturbs the
conductivity $\sigma (x)$. Following the observations made in \cite{AE2,AE3}%
, we assume that the perturbation is proportional to the local value of the
conductivity; thus, the perturbed conductivity $\sigma ^{new}(x)$ equals to $%
\sigma (x)\exp (\eta (x))$, where the perturbation exponent $\eta(x)$ is such that
$|\eta (x)|\ll 1$ and is compactly supported. Let $u^{new}(x)=u(x)+w_{\eta }(x)$
be the potential corresponding to the perturbed conductivity $\sigma ^{new}(x)$
and $w_{\eta }(x)$ be the perturbation thereof. By substituting these perturbed
values into (\ref{original1}) one obtains
\begin{equation}
\Delta \left[ u(x)+w_{\eta }(x)\right] +\nabla \left[ u(x)+w_{\eta }(x)%
\right] \cdot \nabla \left[ \ln \sigma (x)+\eta (x)\right] =0.
\label{original2}
\end{equation}
Further, by neglecting second order terms (in $\eta $) and by subtracting (%
\ref{original1})\ from (\ref{original2}) we arrive at the the following
equation:
\begin{equation}
\Delta w_{\eta }(x)+\nabla w_{\eta }(x)\cdot \nabla \ln \sigma (x)=-\nabla
u(x)\cdot \nabla \eta (x).  \label{original3}
\end{equation}%
Finally, by multiplying (\ref{original3}) by $\sigma (x)$ we find that $%
w_{\eta }(x)$ satisfies equation
\begin{equation}
\label{E:perturb}\nabla\cdot\sigma(x)\nabla w_{\eta}(x)=-\sigma(x)\nabla
u(x)\cdot\nabla\eta(x)
\end{equation}
subject to the homogeneous Neumann boundary conditions. Since the values of
$u(x)$ and $u^{new}(x)$ are measured on the boundary, the Dirichlet data for
$w_{\eta}(x)$ are known. It will be sufficient for our purposes to measure a
certain functional of the boundary values of $w_{\eta}(x)$. Let us fix a
function $I(z)$ defined on $\partial\Omega$, and define the corresponding
measurement functional $M_{I}(\eta)$ as follows:
\begin{equation}
\label{E:functional}M_{I,J}(\eta):= \int\limits_{\partial\Omega}w_{\eta
}(z)I(z)dz.
\end{equation}
Here the subscript $J$ on the left reminds about the dependence of $w$ on the
current~$J$. Function $I(z)$ does not have to be a function in the classical sense;
it may also
be chosen to be a distribution, for example a sum of delta-functions. In the latter
case it would model measurements obtained
by a set of point-like electrodes. Since the data corresponding to all electrodes
then would be added together, the noise sensitivity  of such a scheme is quite low,
and our numerical experiments (not presented here) confirm that.

Our goal is to reconstruct $\sigma(x)$ from measurements of $M_{I,J}(\eta)$
corresponding to a sufficiently rich set of perturbations $\eta(x)$ in
(\ref{E:perturb}).

The simplest case is when one can achieve perfect focusing, and thus $\eta
_{y}(x)\approx C\delta(x-y)$, where the point $y$ scans through $\Omega$. Then
the reconstruction needs to be done from the values
\begin{equation}
M_{I,J,\delta}(y):=\int\limits_{\partial\Omega}w_{\eta_{y},J}(z)I(z)dz.
\notag \end{equation}
However, this assumption of perfect focusing is unrealistic \cite{Localized}.
More realistic are, for instance, mono-chromatic planar or spherical waves, or
spreading spherical fronts. We assume here that ideal point-like transducers
are excited by an infinitesimally short electrical pulse. If we assume (without loss
of generality) that the speed of sound equals~1,
the acoustic
pressure $W_{t,z}(x)$ generated by a transducer placed at point $z$ (outside
$\Omega$) solves the following initial value problem for the wave equation:
\begin{equation}
\begin{cases}
\Delta_x W_{t,z}(x)=\frac{\partial^{2}}{\partial
t^{2}}W_{t,z}(x),\quad x\in\mathbb{R}^{3},\quad t\in\lbrack0,\infty)\\
W_{0,z}(x)=\delta(|x-z|),\\
\frac{\partial}{\partial t}W_{0,z}(x)=0.
\end{cases}
.  \nonumber
\end{equation}
Solution of this problem is well-known \cite{Vladimirov}:
\begin{equation}
W_{t,z}(x)=\frac{\partial}{\partial t}\left(  \frac{\delta(t-|x-z|)}{4\pi
t}\right)  ;\label{E:pulse}%
\end{equation}
it has the form of the propagating spherical front with the radius $t$ centered at
$z$. (The time derivative of the $\delta$-function in (\ref{E:pulse}) results
naturally from the $\delta$-excitation of the transducer; the spherical waves
we used in \cite{KuKuAET} can be obtained by anti-differentiation of the
signal corresponding to (\ref{E:pulse}).)\footnote{Other \textquotedblleft
bases\textquotedblright\ of waves, e.g. radial mono-chromatic, or planar could
also be used \cite{KuKuAET}.}

The perturbation $\eta_{t,z}(x)$ of the conductivity caused by the propagating
front $W_{t,z}(x)$ equals $\eta_{0}W_{t,z}(x)$, where $\eta_{0}$ is some small
fixed proportionality constant (reflecting the smallness of the
acousto-electric effect). The corresponding measurements then are (after
factoring out $\eta_{0}$):
\begin{equation}
\label{E:spherical}M_{I,J}(t,z):=\int\limits_{\partial\Omega}w_{W_{t,z}
,J}(z)I(z)dz.
\end{equation}
Due to the linear dependence of the measurements on the acoustic perturbation
$\eta$, one can try to do a ``basis change'' type of calculation, which would
produce the ``focused'' data $M_{I,J,\delta}(y)$ from the more realistic
``non-focused'' measurements $M_{I,J}(t,z)$. In particular, as it is explained
in \cite{KuKuAET,Oberw}, if one knows the data (\ref{E:spherical}) for all
$t\in [0,\infty]$ and $z\in\Sigma$ (where $\Sigma$ is a closed curve (surface
in $3D$) surrounding $\Omega)$, then $M_{I,J,\delta}(y)$ can be reconstructed
by methods of thermoacoustic tomography. In particular, if $\Sigma$ is a
sphere, circle, cylinder, or a surface of a cube, explicit inversion formulas
exist that can recover $M_{I,J,\delta}(y)$ (see \cite{KuKu}). For general
closed surfaces, other efficient methods exist (e.g., time reversal). This
transformation is known to be stable. In fact, as it will be explained below, in
the version of synthetic focusing used here, it is smoothing.

We thus assume that $M_{I,J,\delta}(y)$ are known for all $y\in\Omega$, (e.g.
they are obtained by synthetic focusing or by direct measurements.) For our
purposes it will be sufficient to use just two functions $I_{1}(z), I_{2}(z)$
as both the current patterns and the weights in the functionals
(\ref{E:functional}). We thus measure or synthesize the following values:
\begin{equation}
M_{i,j}(y):= \int\limits_{\partial\Omega}w_{\eta_{y},I_{i}}(z)I_{j}
(z)dz,\qquad i,j=1,2. \label{E:delta}
\end{equation}

We now interpret this data in a different manner. Namely, let $u_{j}(x)$,
$j=1,2$ be the solutions of (\ref{original}) corresponding to the boundary
currents (i.e., Neumann data) $I_{j}$. Then
\begin{equation}
\label{E:div}\nabla\cdot\sigma(x)\nabla w_{j, \delta_{y}}(x)=-\sigma(y)\nabla
u_{j}(y)\cdot\nabla\delta(x-y).
\end{equation}
Since
\begin{equation}
M_{i,j}=\int\limits_{\partial\Omega}w_{i}(z)I_{j}(z)dz,
\notag \end{equation}
equation (\ref{E:div}) and the divergence theorem lead to the formula:
\begin{equation}
\label{E:mij}M_{i,j}(x_{0})=\sigma(x_{0})\nabla u_{i}(x_{0})\cdot\nabla
u_{j}(x_{0}).
\end{equation}

Thus, for any interior point $x\in\Omega$ and any two current profiles $I_{j},
j=1,2$ on the boundary, the values of the expressions (\ref{E:mij}) can be
extracted from the measured data

Our goal now is to try to recover the conductivity from these values. The same
problem in $2D$ was addressed in \cite{Cap}, but our approach to
reconstruction is different.

\section{Reconstructing the $2D$ conductivity from focused data using two
currents} \label{S:reconstruction}

We will assume here availability of the measurement data $M_{I,J}(x)$ for all
$x \in\Omega$, no matter whether they were obtained by applying focused beams,
or by synthetic focusing. We will consider now the situation where the
conductivity $\sigma(x)$ is considered to be a (relatively) small perturbation
of a known benchmark conductivity $\sigma_{0}(x)$:
\begin{equation}
\sigma(x)=\sigma_{0}(x)(1+\varepsilon\rho(x)),
\end{equation}
where $\varepsilon\ll1$ and $\rho=0$ near the boundary of the domain.
(Numerical experiments show that our method yields quite accurate
reconstructions even when the true conductivity differs significantly from the
initial guess $\sigma_{0}$).

It will be also assumed that two distinct current patterns $I_{j}$, $j=1,2$ on
the boundary are fixed, and the two resulting potentials $u_{j}, j=1,2$ with
the benchmark conductivity $\sigma_{0}$:
\begin{equation}
\nabla\cdot\sigma_{0}(x)\nabla u_{j}(x)=0
\notag \end{equation}
corresponding to the two prescribed sets of boundary currents. These
potentials can be computed and are assumed to be known.

Correspondingly, the unknown true potentials $w_{j}(x)=u_{j}(x)+\varepsilon
v_{j}(x)+o(\varepsilon)$ for the actual conductivity $\sigma$ satisfy the
equations
\begin{equation}
\nabla\cdot\sigma\nabla(u_{j}+\varepsilon v_{j})=0
\notag \end{equation}
with the same boundary currents as $u_{j}$.

According to the discussion in the previous section, using acoustic
delta-perturbations (real or synthesized), we can obtain for any point $x$ in
the domain $\Omega$ the values
\begin{equation}
\label{E:M0}M^{0}_{j,k}(x):=\sigma_{0}(x)\nabla u_{j}(x)\cdot\nabla u_{k}(x),
\end{equation}
which can be computed numerically using the background conductivity
$\sigma_{0}$, and
\begin{equation}
\label{E:M}M_{j,k}(x):=\sigma(x)\nabla w_{j}(x)\cdot\nabla w_{k}
(x)=M^{0}_{j,k}+\varepsilon g_{j,k}+o(\varepsilon),
\end{equation}
which are obtained by boundary measurements. Now we can forget about the
acoustic modulation and concentrate on reconstructing $\rho(x)$ (and thus
$\sigma(x)$) from the known $M_{j,k}(x)$, or, neglecting higher order terms,
from $g_{j,k}(x)$.

Let us re-write (\ref{E:M}) in the following form:%
\begin{equation}
\sigma (x)\nabla \left[ u_{j}(x)+\varepsilon v_{j}(x)\right] \cdot \nabla %
\left[ u_{k}(x)+\varepsilon v_{k}(x)\right] =M_{j,k}^{0}+\varepsilon
g_{j,k}+o(\varepsilon ),  \label{E:M1}
\end{equation}%
By subtracting (\ref{E:M0}) from (\ref{E:M1}) one obtains formulas
\begin{equation}
g_{j,k}(x)=\sigma\left(  \nabla u_{j}\cdot\nabla v_{k}+\nabla u_{k}\cdot\nabla
v_{j}\right)  +o(\varepsilon).\label{E:gjk}
\end{equation}
We will drop the $o(\varepsilon)$ terms in the following calculations. We
introduce the new vector fields $U_{j}=\sqrt{\sigma_{0}}\nabla u_{j}$ and
$W_{j}=\sqrt{\sigma}\nabla(u_{j}+\varepsilon v_{j})=U_{j}+\varepsilon V_{j},$
so that
\begin{equation}
\nabla\cdot\sqrt{\sigma_{0}}U_{j}=0
\notag \end{equation}
and
\begin{equation}
\nabla\cdot\sqrt{\sigma}W_{j}=0.
\notag \end{equation}
We would like to find $W_{j}$. The last equation can be re-written, taking
into account that, up to $o(\varepsilon)$ terms, $\sqrt{\sigma}\approx
\sqrt{\sigma_{0}}(1+\frac{1}{2}\varepsilon\rho)$ and $\ln\sigma=\ln\sigma
_{0}+\varepsilon\rho$, as follows:
\begin{equation}
\nabla\cdot\sqrt{\sigma_{0}}(1+\varepsilon\rho/2)(U_{j}+\varepsilon V_{j})=0
\notag \end{equation}
or
\begin{equation}
\nabla\cdot(U_{j}+\varepsilon V_{j})+\frac{1}{2}(U_{j}+\varepsilon V_{j}
)\cdot\nabla(\ln\sigma+\varepsilon\rho)=0.
\notag \end{equation}
By collecting the terms of the zero and first order in $\varepsilon$ we obtain
\begin{equation}
\nabla\cdot U_{j}+\frac{1}{2}U_{j}\cdot\nabla\ln\sigma=0
\notag \end{equation}
and
\begin{equation}
\nabla\cdot V_{j}+\frac{1}{2}U_{j}\cdot\nabla\rho+\frac{1}{2}V_{j}\cdot
\nabla\ln\sigma=0
\notag \end{equation}
or
\begin{equation}
\nabla\cdot V_{j}+\frac{1}{2}V_{j}\cdot\nabla\ln\sigma=-\frac{1}{2}U_{j}
\cdot\nabla\rho.
\notag \end{equation}
Equivalently
\begin{equation}
\nabla\cdot\sqrt{\sigma}V_{j}=-\frac{1}{2}\sqrt{\sigma}U_{j}\cdot\nabla\rho.
\notag \end{equation}
With this new notation, the measurements can be expressed (neglecting higher
order terms) as follows:
\begin{equation}
(U_{j}+\varepsilon V_{j})\cdot(U_{k}+\varepsilon V_{k})=M_{j,k}=M^{0}
_{j,k}+\varepsilon g_{j,k},
\notag \end{equation}
which leads to
\begin{align*}
U_{j}\cdot U_{k}  &  =M^{0}_{j,k},\\
U_{j}\cdot V_{k}+U_{k}\cdot V_{j}  &  =g_{j,k}.
\end{align*}
In particular, we arrive to three independent equations for $V_{j}$:
\begin{align}
U_{1}\cdot V_{1}  &  =g_{1,1}/2\nonumber\\
U_{2}\cdot V_{2}  &  =g_{2,2}/2\label{E:mainsys}\\
U_{1}\cdot V_{2}+U_{2}\cdot V_{1}  &  =g_{1,2}.\nonumber
\end{align}

These equations will be our starting point for deriving reconstruction
algorithms, as well as uniqueness and stability results.

We consider now the case when the benchmark conductivity (initial conductivity
guess) is constant: $\sigma_{0}(x)\equiv1$.

\subsection{The constant benchmark conductivity $\sigma_{0} (x)=1$}
\label{S:constant}

We will choose the boundary currents $\frac{\partial}{\partial n}u_{j}(x)$ to be
equal to $n(x)\cdot e_{j}$, where $n(x)$ is the unit external normal to the
boundary and $e_{1}=(1,0),e_{2}=(0,1)$ are the canonical basis vectors. Then for
the conductivity $\sigma_{0}=1$ the resulting potentials $u_{j}(x)$ are equal to
$x_{j},$ and the fields $U_{j}$ are equal to $e_{j}$:
\begin{equation}
U_{j}=\nabla u_{j}=e_{j},j=1,2.\nonumber
\end{equation}

We thus obtain formulas
\begin{equation}
\label{E:Eq_gjk}\left\{
\begin{array}
[c]{c}
2\frac{\partial v_{1}}{\partial x_{1}}+\rho=g_{1,1}\\
2\frac{\partial v_{2}}{\partial x_{2}}+\rho=g_{2,2}\\
\frac{\partial v_{1}}{\partial x_{2}}+\frac{\partial v_{2}}{\partial x_{1}
}=g_{1,2}
\end{array}
\right.
\end{equation}
as well as the equations
\begin{equation}
\Delta v_{j}=-\frac{\partial}{\partial x_{j}}\rho,\qquad j=1,2.
\label{E:laplacians}
\end{equation}
Differentiating the equations (\ref{E:Eq_gjk}), we obtain
\begin{equation}
\left\{
\begin{array}
[c]{c}
2\frac{\partial^{2}v_{1}}{\partial x_{1}^{2}}+\frac{\partial}{\partial x_{1}
}\rho=\frac{\partial}{\partial x_{1}}g_{1,1}\\
2\frac{\partial^{2}v_{1}}{\partial x_{1}\partial x_{2}}+\frac{\partial
}{\partial x_{2}}\rho=\frac{\partial}{\partial x_{2}}g_{1,1}\\
2\frac{\partial^{2}v_{2}}{\partial x_{2}^{2}}+\frac{\partial}{\partial x_{2}
}\rho=\frac{\partial}{\partial x_{2}}g_{2,2}\\
2\frac{\partial^{2}v_{2}}{\partial x_{1}\partial x_{2}}+\frac{\partial
}{\partial x_{1}}\rho=\frac{\partial}{\partial x_{1}}g_{2,2}\\
\frac{\partial^{2}v_{1}}{\partial x_{1}\partial x_{2}}+\frac{\partial^{2}
v_{2}}{\partial x_{1}^{2}}=\frac{\partial}{\partial x_{1}}g_{1,2}\\
\frac{\partial^{2}v_{1}}{\partial x_{2}^{2}}+\frac{\partial^{2}v_{2}}{\partial
x_{1}\partial x_{2}}=\frac{\partial}{\partial x_{2}}g_{1,2}
\end{array}
\right. \label{E:bigsys}
\end{equation}
Combining the 2nd, 3rd, and 5th equations in (\ref{E:bigsys}), we arrive to
\begin{equation}
0=\frac{\partial}{\partial x_{2}}g_{1,1}-2\frac{\partial}{\partial x_{1}
}g_{1,2}-\frac{\partial}{\partial x_{2}}g_{2,2}+2\Delta v_{2}.
\notag \end{equation}
Utilizing (\ref{E:laplacians}) with $j=2$ and differentiating with respect to
$x_{2}$, we obtain
\begin{equation}
\frac{\partial^{2}}{\partial x_{2}^{2}}\rho=\frac{1}{2}\frac{\partial^{2}
}{\partial x_{2}^{2}}(g_{1,1}-g_{2,2})-\frac{\partial^{2}}{\partial
x_{1}\partial x_{2}}g_{1,2}.
\notag \end{equation}
Similarly,
\begin{equation}
\frac{\partial^{2}}{\partial x_{1}^{2}}\rho=\frac{1}{2}\frac{\partial^{2}
}{\partial x_{1}^{2}}(g_{2,2}-g_{1,1})-\frac{\partial^{2}}{\partial
x_{1}\partial x_{2}}g_{1,2}.
\notag \end{equation}
Adding the last two equalities, we obtain the Poisson type equation
\begin{equation}
\Delta\rho=\frac{1}{2}\left(  \frac{\partial^{2}}{\partial x_{1}^{2}}
-\frac{\partial^{2}}{\partial x_{2}^{2}}\right)  (g_{2,2}-g_{1,1}
)-2\frac{\partial^{2}}{\partial x_{1}\partial x_{2}}g_{1,2}
\label{E:smartPoisson}
\end{equation}
for the unknown function $\rho$. Notice that all expressions in the right hand
side are obtained from the measured data and that by our assumption $\rho$
satisfies the zero Dirichlet condition at the boundary.

This reduction clearly allows for algorithmic reconstruction, as well as
proving (under appropriate smoothness assumptions on $\sigma$) local
uniqueness and Lipschitz stability of reconstruction (see
Section~\ref{S:stability}).

\subsection{A parametrix solution for smooth benchmark conductivity
$\sigma_{0}(x)$}\label{S:smooth}

We would like to present now a sometimes useful observation for the situation when
benchmark conductivity $\sigma_{0}$ is smooth, but not necessarily constant
(e.g., a standard EIT reconstruction would provide such an approximation). In
this case, we will find a parametrix solution, i.e. will determine $\sigma(x)$
up to smoother terms.

As it has already been discussed, perturbation $\varepsilon v_{j}$ of the
potential $u_{j}$ satisfies the equation
\begin{equation}
\nabla\cdot\sigma_{0}\nabla v_{j}=-\sigma_{0}\nabla u_{j}\cdot\nabla\rho.
\notag \end{equation}
Since $\sigma$ is smooth and non-vanishing, up to smoother terms we can write
\begin{equation}
\Delta v_{j}\thickapprox-\nabla u_{j}\cdot\nabla\rho
\notag \end{equation}
and
\begin{equation}
v_{j}\thickapprox-(\nabla u_{j}\cdot\nabla)(\Delta^{-1}\rho)
\notag \end{equation}
where $\Delta^{-1}$ is the inverse to the Dirichlet Laplacian in $\Omega$.
Again up to smoother terms, we have
\begin{align*}
U_{k}\cdot V_{j}  &  =\sqrt{\sigma}\nabla u_{k}\cdot\sqrt{\sigma}(\rho/2\nabla
u_{j}+\nabla v_{j})\\
&  =\sigma\rho/2\nabla u_{k}\cdot\nabla u_{j}+\sigma(\nabla u_{k}\cdot
\nabla)(\nabla u_{j}\cdot\nabla)\Delta^{-1}\rho.
\end{align*}
The latter expression is symmetric up to smoothing terms and equations
(\ref{E:mainsys}) can be re-written as
\begin{align*}
U_{1}\cdot V_{1}  &  =g_{1,1}/2\\
U_{2}\cdot V_{2}  &  =g_{2,2}/2\\
U_{1}\cdot V_{2}  &  =g_{1,2}/2+\mbox{ a smoother term}\\
U_{2}\cdot V_{1}  &  =g_{1,2}/2+\mbox{ a smoother term}.
\end{align*}
Under such an approximation, assuming that currents $\nabla u_{1}$ and $\nabla
u_{2}$ are not parallel, which is known to be possible to achieve \cite{Ales},
one can recover $\varepsilon V_{j}$ at each point $x$. Therefore, (more)
accurate solutions $W_{j}=U_{j}+\varepsilon V_{j}$ can be found. We note that
$\nabla\cdot\sqrt{\sigma}W_{j}=0 $ and so
\begin{equation}
W_{j}\cdot\nabla\ln\sigma=-2\nabla\cdot W_{j}.
\notag \end{equation}
On the other hand, since $W_{j}=\sqrt{\sigma}\nabla(u_{j}+\varepsilon v_{j}),$
we have
\begin{equation}
\nabla\times\frac{W_{j}}{\sqrt{\sigma}}=0.
\notag \end{equation}
This can be re-written as
\begin{equation}
W_{j}\times\nabla\ln\sigma=-2\nabla\times W_{j}
\notag \end{equation}
or
\begin{equation}
W_{j}^{\perp}\cdot\nabla\ln\sigma=-2\nabla\times W_{j},
\notag \end{equation}
where $W_{j}^{\perp}$ is the vector obtained from $W_{j}$ by the
counter-clockwise $90^{o}$ rotation (i.e. $W_{j}^{\perp}\cdot W_{j}=0$ and
$|W_{j}^{\perp}|=|W_{j}|)$.

Since for each $j=1,2$ vectors $W_{j}$ and $W_{j}^{\perp}$ form an orthogonal
basis, one has
\begin{equation}
\nabla\ln\sigma=-\frac{2}{|W_{j}|^{2}}(W_{j}^{\perp}(\nabla\times W_{j}
)+W_{j}(\nabla\cdot W_{j})),
\notag \end{equation}
and thus
\begin{equation}
\Delta\ln\sigma=-\operatorname{div}\frac{2}{|W_{j}|^{2}}(W_{j}^{\perp}
(\nabla\times W_{j})+W_{j}(\nabla\cdot W_{j})).
\notag \end{equation}
We compute now $\ln\sigma$ by taking the average of the two values of $j$ and
then solving the Poisson equation
\begin{equation}
\Delta\ln\sigma=-\operatorname{div}\sum\limits_{j=1}^{2}\frac{2}{|W_{j}|^{2}
}(W_{j}^{\perp}(\nabla\times W_{j})+W_{j}(\nabla\cdot W_{j})).
\notag \end{equation}

It is interesting to note that this solution reduces to (\ref{E:smartPoisson})
when $\sigma=1$, although (\ref{E:smartPoisson}) holds exactly, not just up to
smoother terms.

\section{Uniqueness and stability}\label{S:stability}
In this section we will assume that $\sigma\in C^{1,\alpha}(\Omega)$, and thus
$\rho$ belongs to this space as well (recall that $\rho$ also vanishes in a fixed
neighborhood of $\partial\Omega$).

The questions of uniqueness and stability in the situation close to ours have
already been addressed in \cite{Cap,bon}, so we will be brief here. Although
considerations of \cite{Cap,bon} were provided in $2D$, the conclusion in our
situation works out the same way in $3D$ if three currents are used.

The standard elliptic regularity \cite{Gilbarg} implies
\begin{proposition}\cite{Cap,bon}
\begin{enumerate}
\item The data $g_{i,j}$ in (\ref{E:M}) determine the conductivity $\sigma=1+\rho$
uniquely.
\item The mappings $\rho(x) \mapsto \{g_{i,j}(x)\}$ of the space
$C^{1,\alpha}_0(\overline{V})$, where $V$ is a compact sub-domain of $\Omega$, are
Fr\'{e}chet differentiable.
\end{enumerate}
\end{proposition}
This justifies our formal linearization near the benchmark  conductivity
$\sigma_0$. Now, the calculations of the Section \ref{S:constant} provide explicit
formulas for the Fr\'{e}chet derivative of the proposition\footnote{In fact, these
formulas easily imply the statement of the proposition in our particular case.}.
In particular,
\begin{equation}
\begin{array}{l}
\frac{\partial}{\partial x_{1}}\rho=\frac{1}{2}\frac{\partial
}{\partial x_{1}}(g_{2,2}-g_{1,1})-\frac{\partial}{\partial x_{2}}g_{1,2},\\
\frac{\partial}{\partial x_{2}}\rho=\frac{1}{2}\frac{\partial
}{\partial x_{2}}(g_{1,1}-g_{2,2})-\frac{\partial}{\partial
x_{1}}g_{1,2}.
\end{array}
\end{equation}
These formulas and vanishing of $\rho$ near $\partial\Omega$ show that the norm of
$\rho$ in $C^{1,\alpha}$ can be estimated from above by such norms of the
functions $\{g_{11},g_{12},g_{22}\}$. In other words, the Fr\'{e}chet derivative
of the mapping
\begin{equation}
\rho\mapsto \{g_{11},g_{12},g_{22}\}
\label{E:mapping}
\end{equation}
is a semi-Fredholm operator with zero kernel. Then the standard implicit function
type argument shows (see, e.g., \cite[Corollary 5.6, Ch. I]{Lang}) that
(\ref{E:mapping}) is an immersion.
This proves local uniqueness and stability for the non-linear problem
(analogous result is obtained in $2D$ in \cite{bon}).

Moreover, since our algorithms start with inverting the Fr\'{e}chet derivative,
this reduces near the constant conductivity the non-linear problem to the one with
an identity plus a contraction operator. This explains why the fixed point
iterations in the following sections converge so nicely.

The $3D$ case with three currents works the same way. Similarly to how it is
done in Section~2.1, for a constant conductivity benchmark $\sigma_0$ one can always
find boundary currents that produce fields $U_j = e_j$, $j=1,2,3$. Then, as explained
in Section~\ref{S:3drec}, one obtains an elliptic system of equations
(see equation (\ref{E:3dsystem})) for reconstructing $\rho(x)$.

\begin{figure}[t]
\begin{center}
\begin{tabular}
[c]{ccc}
\includegraphics[width=1.7in,height=1.7in]{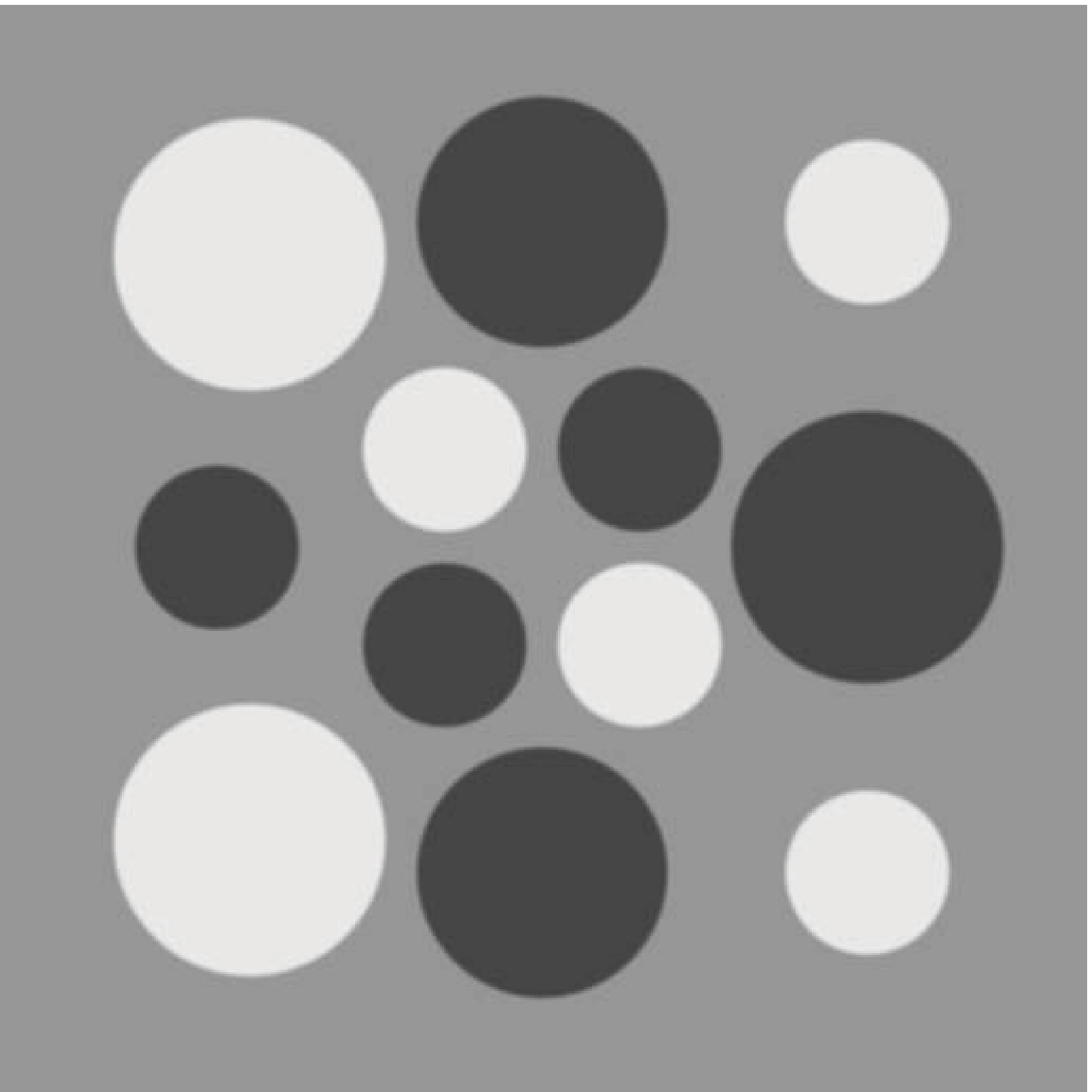} &
\includegraphics[width=1.7in,height=1.7in]{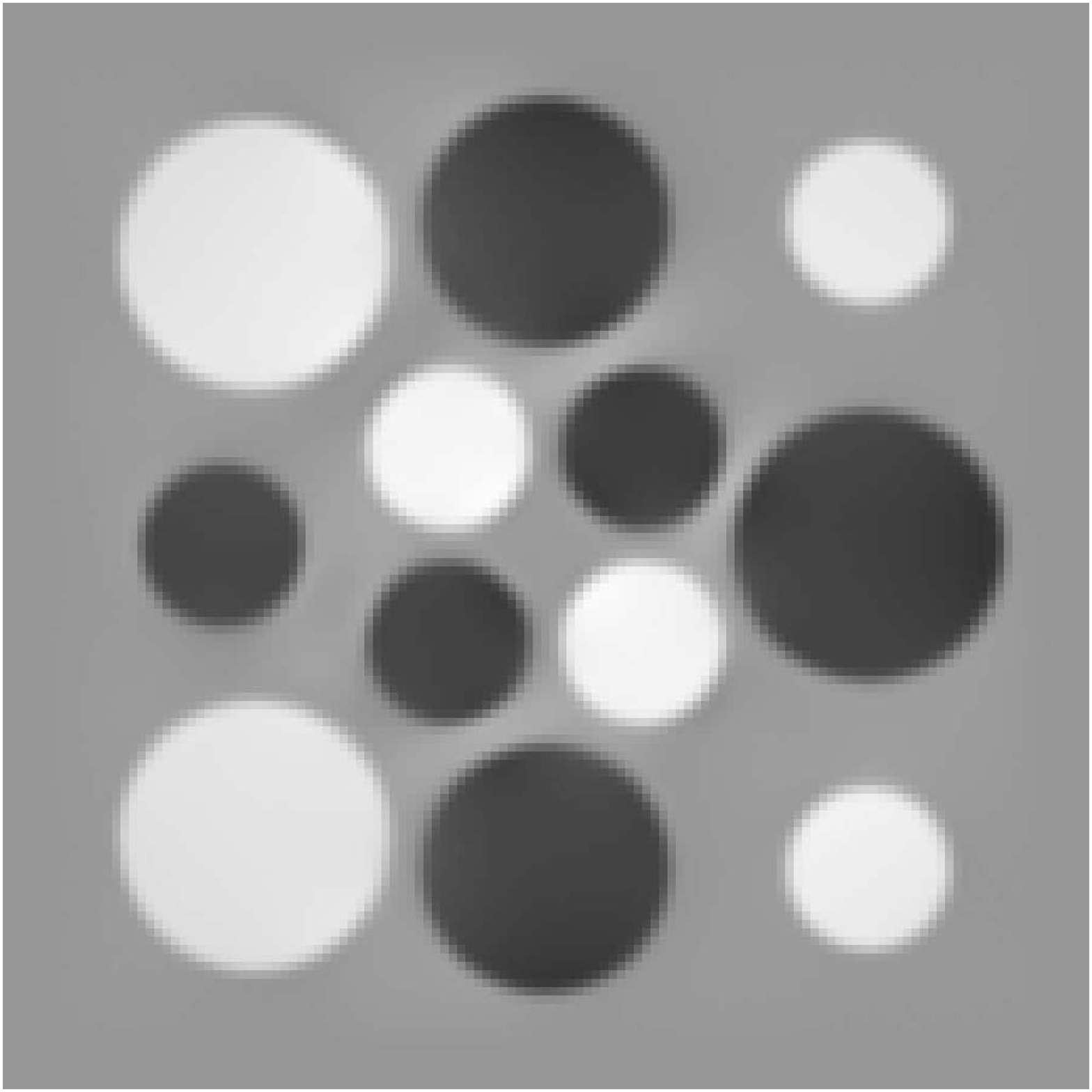} &
\includegraphics[width=1.7in,height=1.7in]{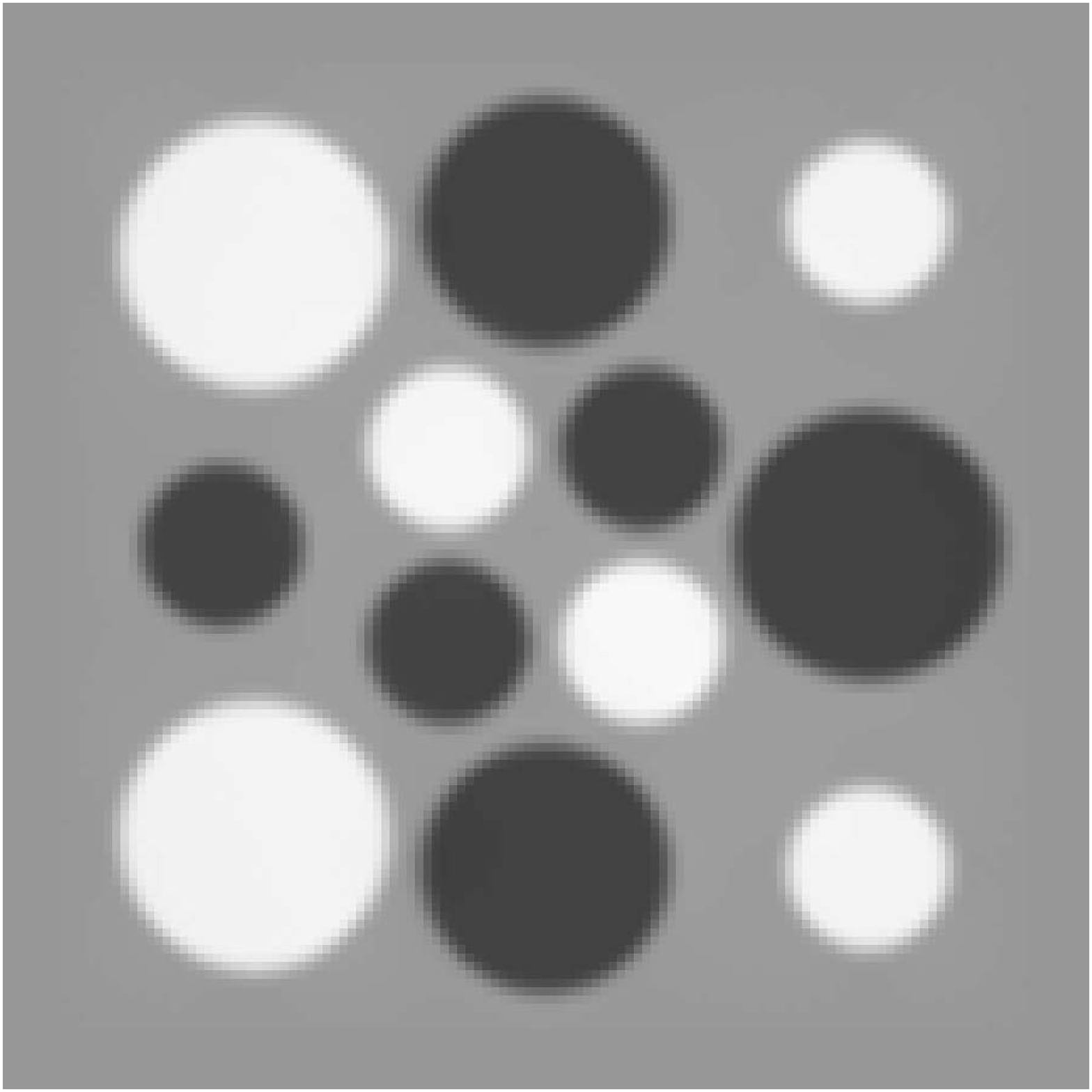}\\
&  & \\
(a) & (b) & (c)
\end{tabular}
\end{center}
\caption{Reconstruction in $2D$ from noiseless data (a)~phantom (b)~iteration~\#0
(c)~iteration~\#1}
\label{F:accucircgray}
\end{figure}

\begin{figure}[t]
\begin{center}
\begin{tabular}
[c]{ccc}
\includegraphics[width=1.7in,height=1.7in]{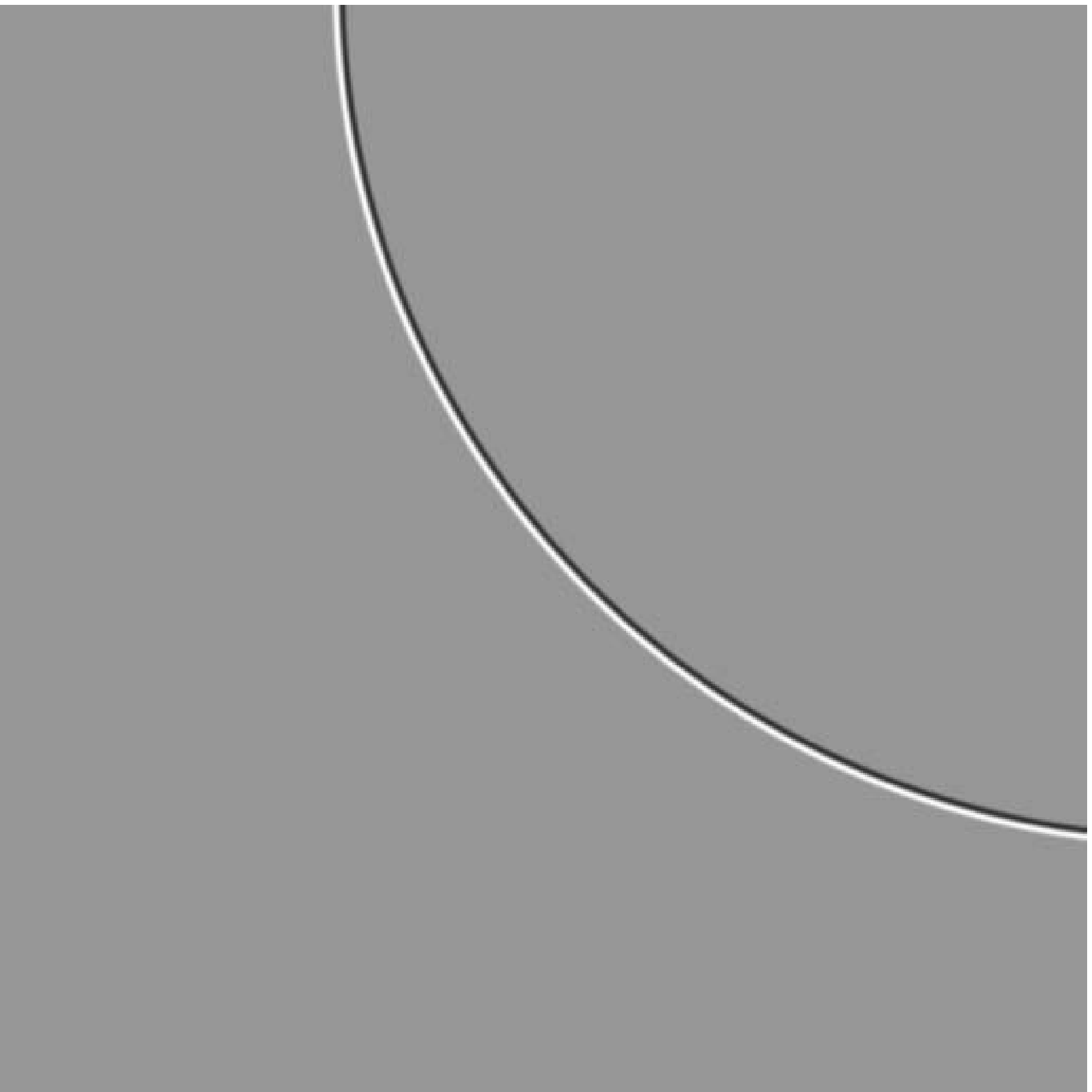} &
\includegraphics[width=1.7in,height=1.7in]{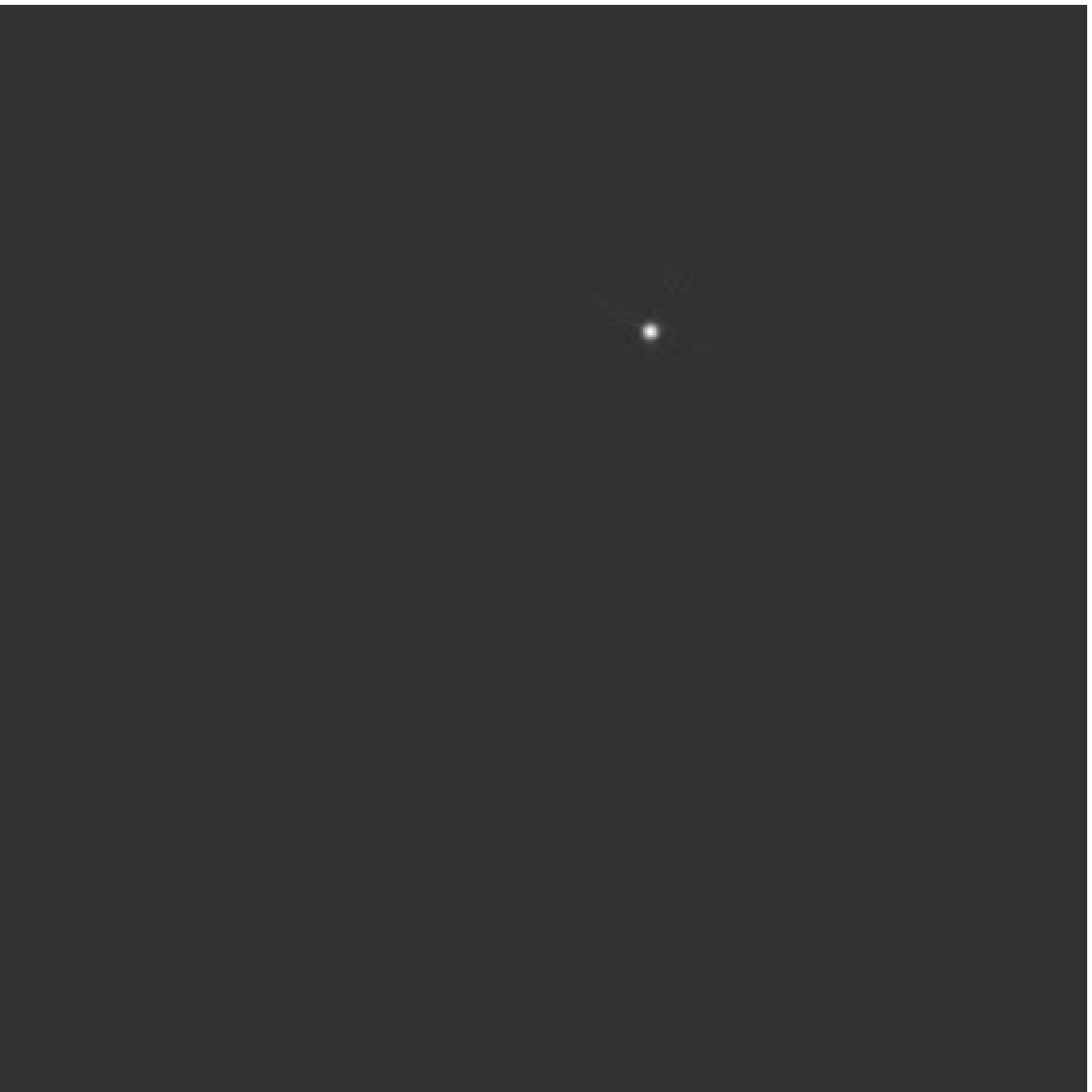} &
\includegraphics[width=1.7in,height=1.7in]{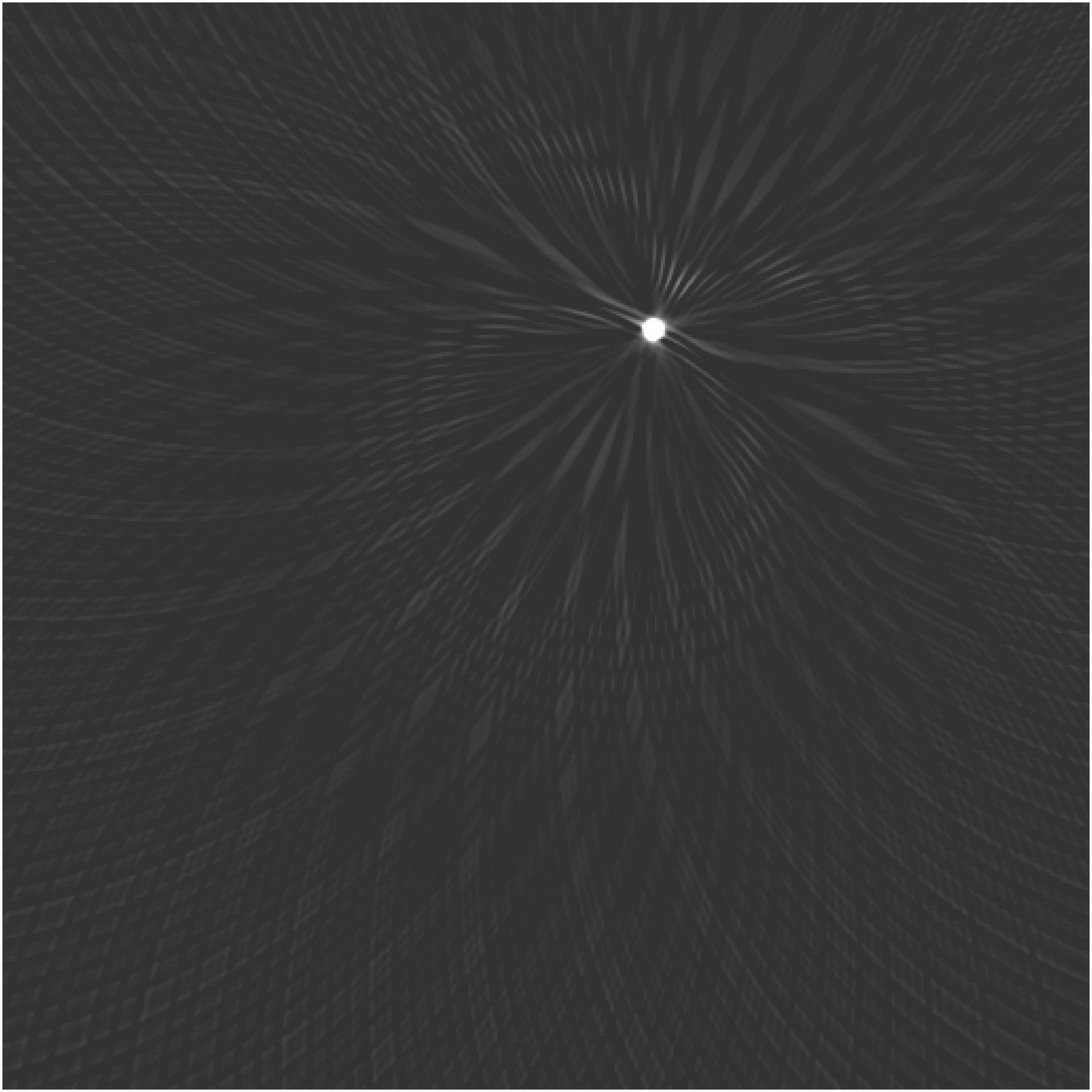}\\
&   \\
(a) & (b) & (c)
\end{tabular}
\end{center}
\caption{(a) Propagating acoustic front (b) the result of focusing at the point
$(0.2,0.4)$ (c) same as (b) with the gray scale showing the lower $10\%$ of the
range of the function}
\label{F:focus}
\end{figure}

\begin{figure}[t]
\begin{center}
\includegraphics[width=2.8in,height=1.10in]{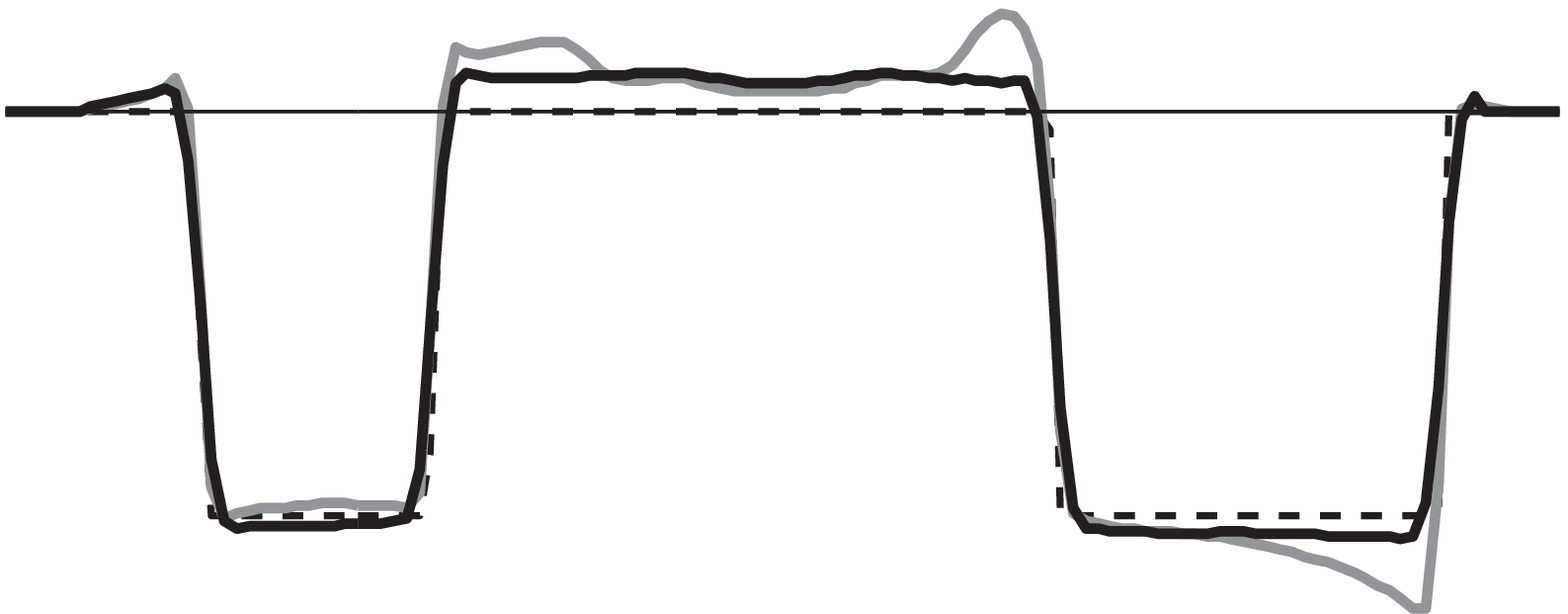}
\end{center}
\caption{Horizontal central cross-section (accurate data): dashed line denotes
the phantom, gray line represents iteration~\# 0, thick black solid line
represents iteration~\#1 }
\label{F:accucircprof}
\end{figure}

\begin{figure}[th]
\begin{center}%
\begin{tabular}
[c]{cc}%
\includegraphics[width=1.7in,height=1.7in]{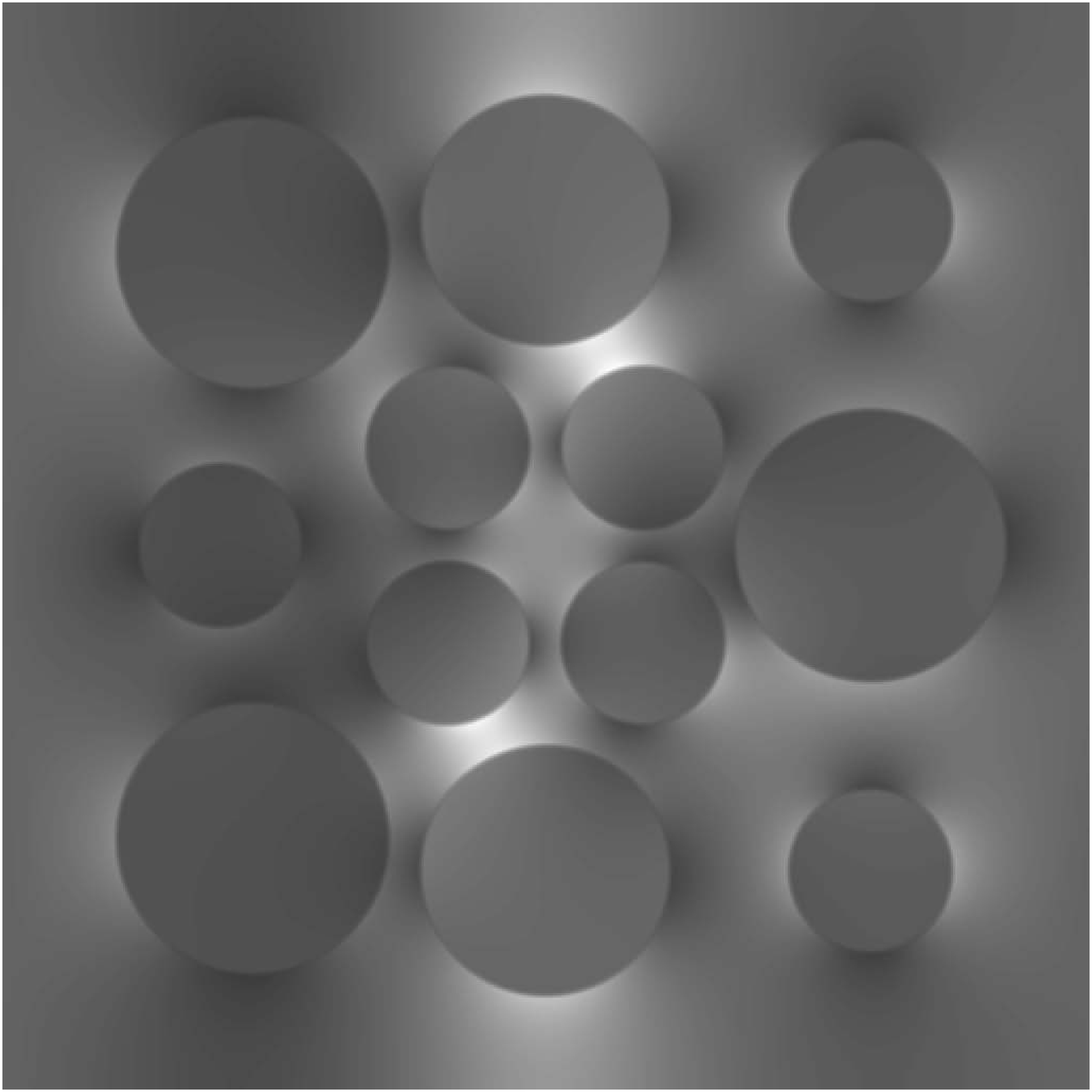} &
\includegraphics[width=1.7in,height=1.7in]{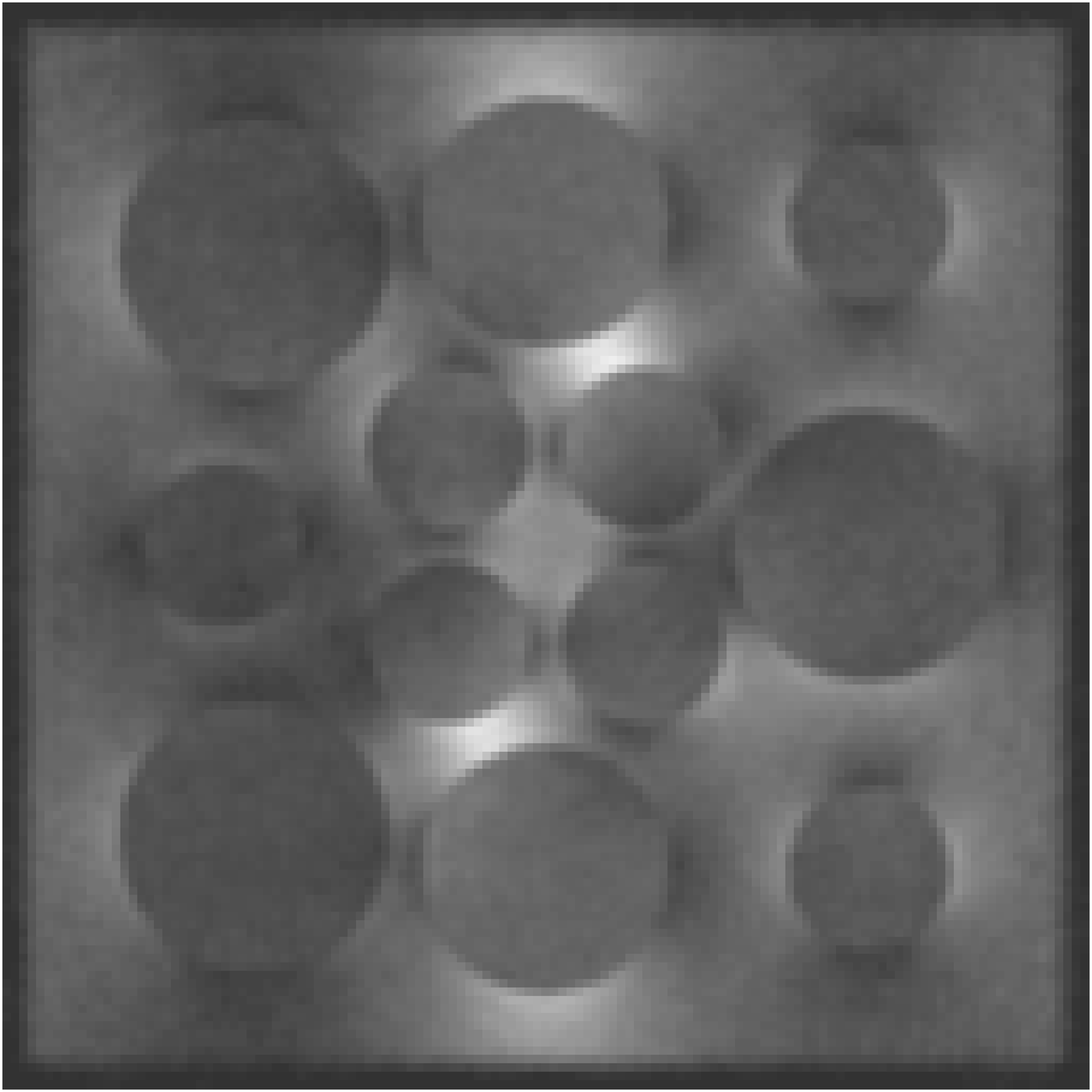}\\

(a) & (b)                                  \\

\includegraphics[width=1.7in,height=1.7in]{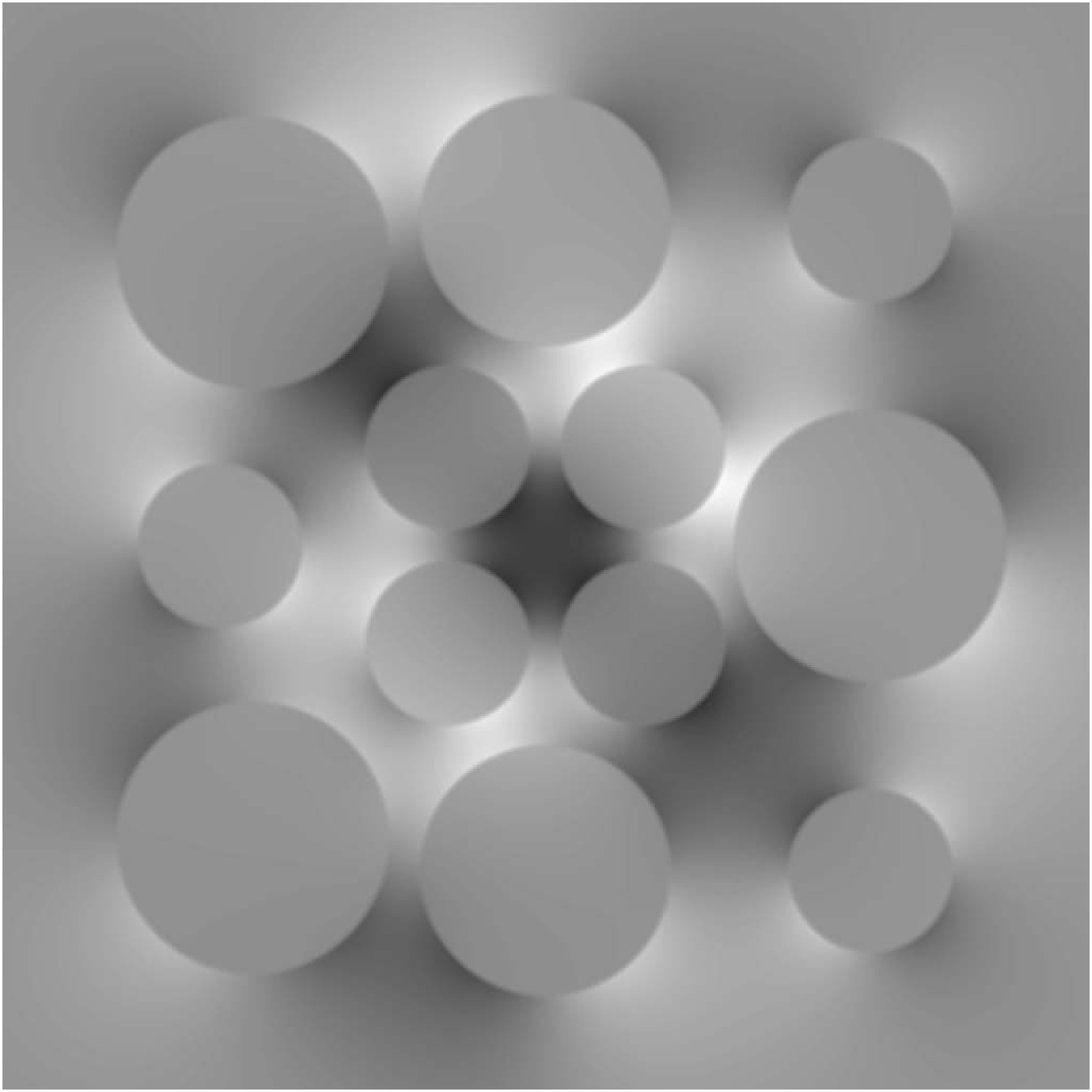} &
\includegraphics[width=1.7in,height=1.7in]{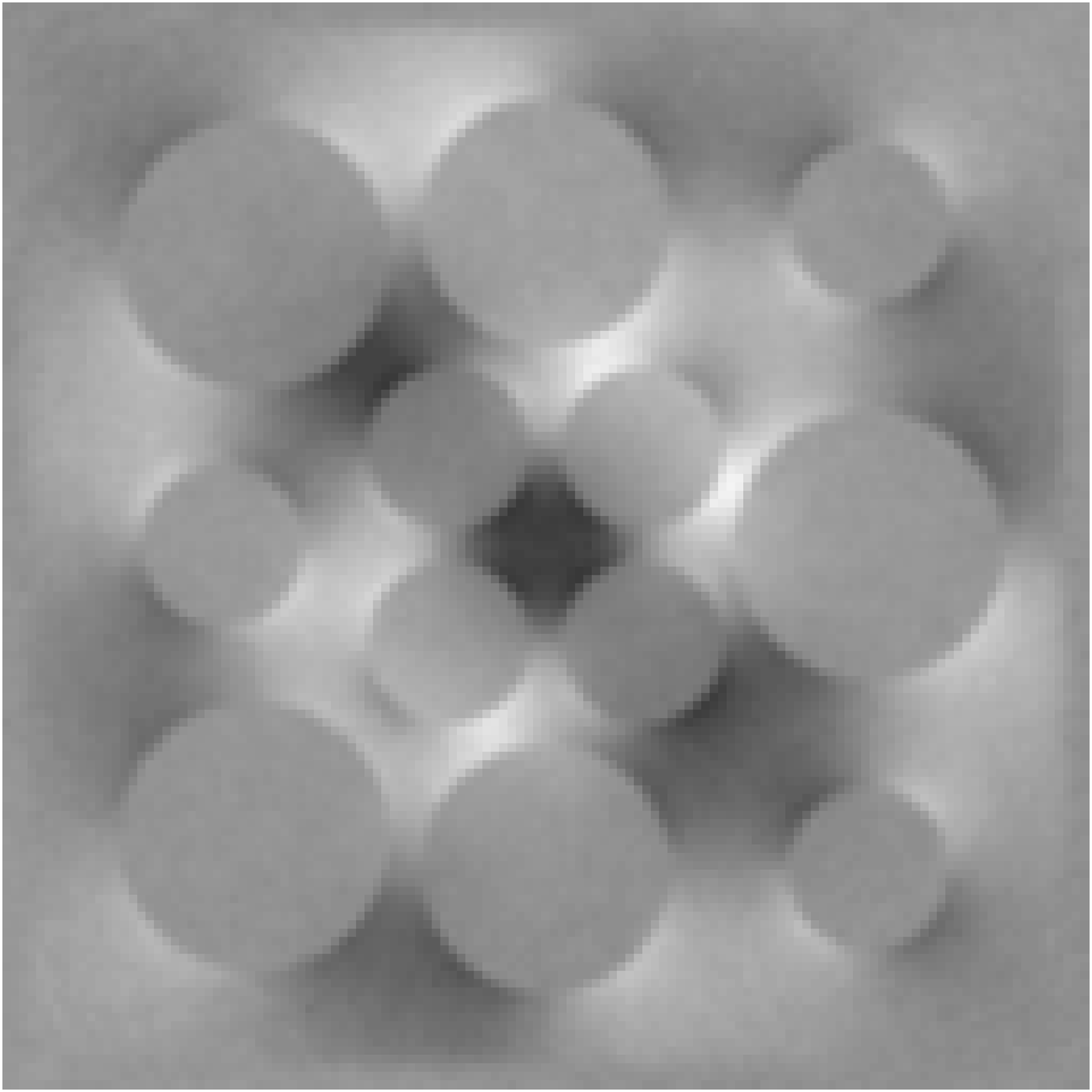}\\

(c) & (d)                         \\
\end{tabular}
\end{center}
\caption{ Functionals $M_{i,j}$: (a)~original $M_{1,1}$ (b)~$M_{1,1}$
reconstructed from data contaminated by $50\%$ noise
(c)~original $M_{1,2}$ (d)~$M_{1,2}$ reconstructed from data
contaminated by $50\%$ noise}
\label{F:mij}
\end{figure}

\section{Numerical examples in $2D$} \label{S:numerical}

We will now illustrate the properties of our algorithm on several
numerical examples in $2D$. Each simulation involves several steps.
First we model the direct problem as follows. For a given phantom of $\sigma$
and a fixed boundary current $J$ we solve equation~(\ref{original})
in the unit square $[-1,1] \times [-1,1]$,
and (for a chosen weight function $I$) we compute the unperturbed boundary
functionals $M^{\mathrm{unperturbed}}_{I,J}$:
\begin{equation}
M^{\mathrm{unperturbed}}_{I,J}:=\int\limits_{\partial\Omega}u(z)I(z)dz.
\end{equation}
Next, for a set of values of $t$ and $z$ we perturb $\sigma$ by multiplying it
by $ \operatorname{exp}(\eta_{t,z}(x))$ with $\eta_{t,z}(x)$ proportional to
the propagating acoustic pulse $W_{t,z}$ given by equation~(\ref{E:pulse}). (In
simulation we used a mollified version of the delta-function, which corresponds
to a transducer with a finite  bandwidth.) For each perturbed $\sigma$ we again
solve equation~(\ref{original}), obtain the solution $u^{\mathrm{perturbed}}$,
and compute functionals
\begin{equation}
M^{\mathrm{perturbed}}_{I,J}(t,z):=\int\limits_{\partial\Omega}
u^{\mathrm{perturbed}}(z)I(z)dz.
\end{equation}
Finally, the difference of $M^{\mathrm{perturbed}}_{I,J}(t,z)$ and
$M^{\mathrm{unperturbed}}_{I,J}$ yields the values of the functionals
$M_{I,J}(t,z)$ given by equation~(\ref{E:spherical}) which we consider the
simulated measurements and the starting point for solving the inverse problems.
In some of our numerical experiments we add values of a random variable to
these functions to simulate the noise in the measurements.

The advantage of computing $M_{I,J}(t,z)$  as the difference of two solutions (as
opposed to obtaining it from the linearized equation~(\ref{E:spherical})) consists
in eliminating  the chance of committing ``an inverse crime''. However, since
subtraction of two numerically computed functions that differ very little
can significantly amplify the relative error, our forward solver has to be very
accurate. In order to achieve high accuracy we approximated the potentials
in the square by Fourier series and used the Fast Fourier transforms (FFT)
to compute the corresponding differential operators. In turn, the application of the
FFTs  allowed us to use fine discretization grids ($513 \times 513$), which,
in combination with smoothing of the simulated $\sigma(x)$ yields the
desired high accuracy. (Such algorithms combining the use of global bases
(such as the trigonometric basis utilized here) with enforcing the equation
in the nodes of the computational grid are called
\emph{pseudospectral}~\cite{Fornberg}; they are very efficient when the
computational domain is simple (e.g. a square) and the coefficients of the
equation are smooth.)

After the measurement data have been simulated, the inverse problem of AET is
solved by reconstructing functions $M_{i,j}$ (see equation~(\ref{E:delta}))
from $M_{I,J} (t,z)$  (synthesis step), and by applying the methods of Section~2
to reconstruct $\varepsilon \rho(x)$ (i.e. the difference
between the true conductivity and the benchmark $\sigma_0$).

\begin{figure}[t]
\begin{center}
\begin{tabular}
[c]{ccc}
\includegraphics[width=1.7in,height=1.7in]{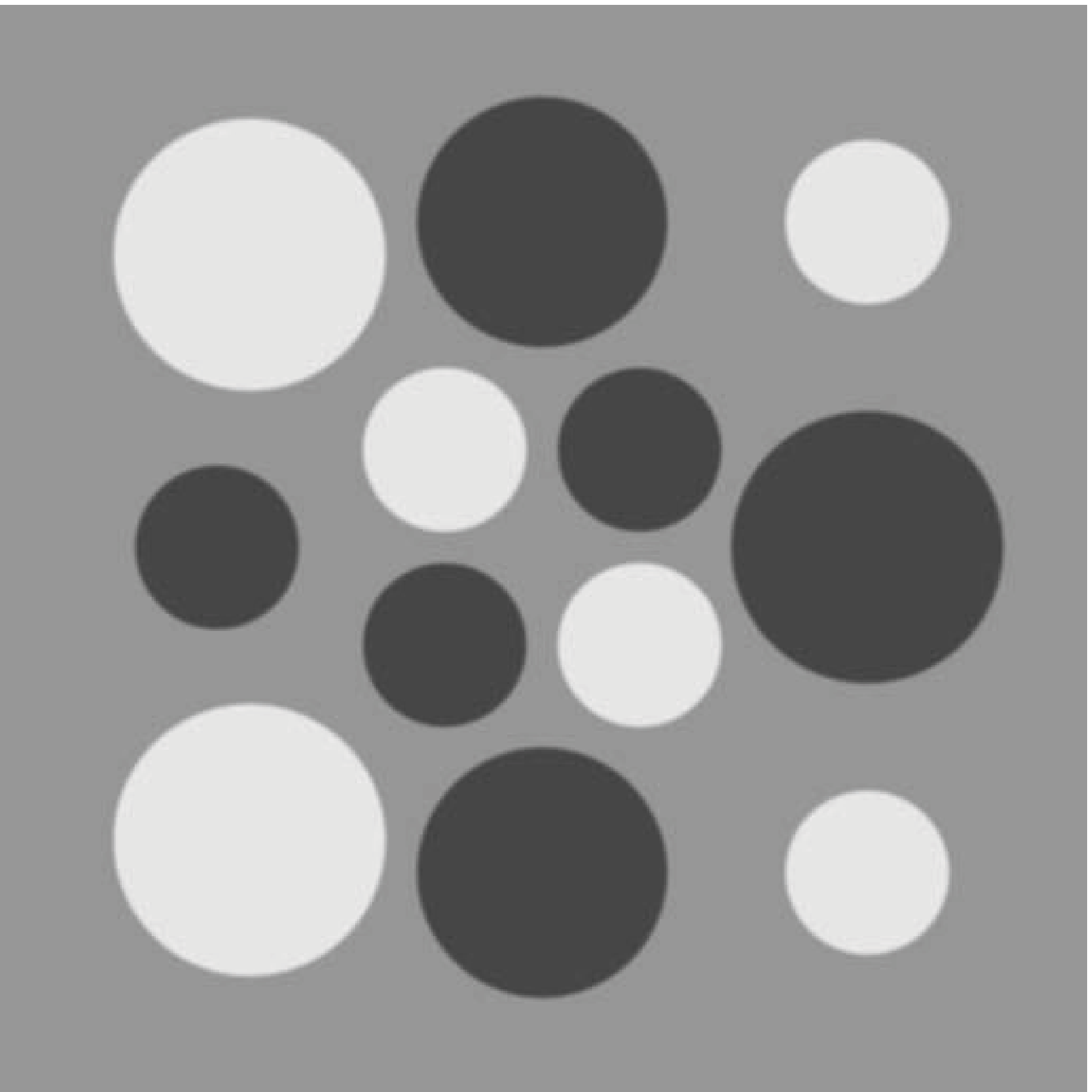} &
\includegraphics[width=1.7in,height=1.7in]{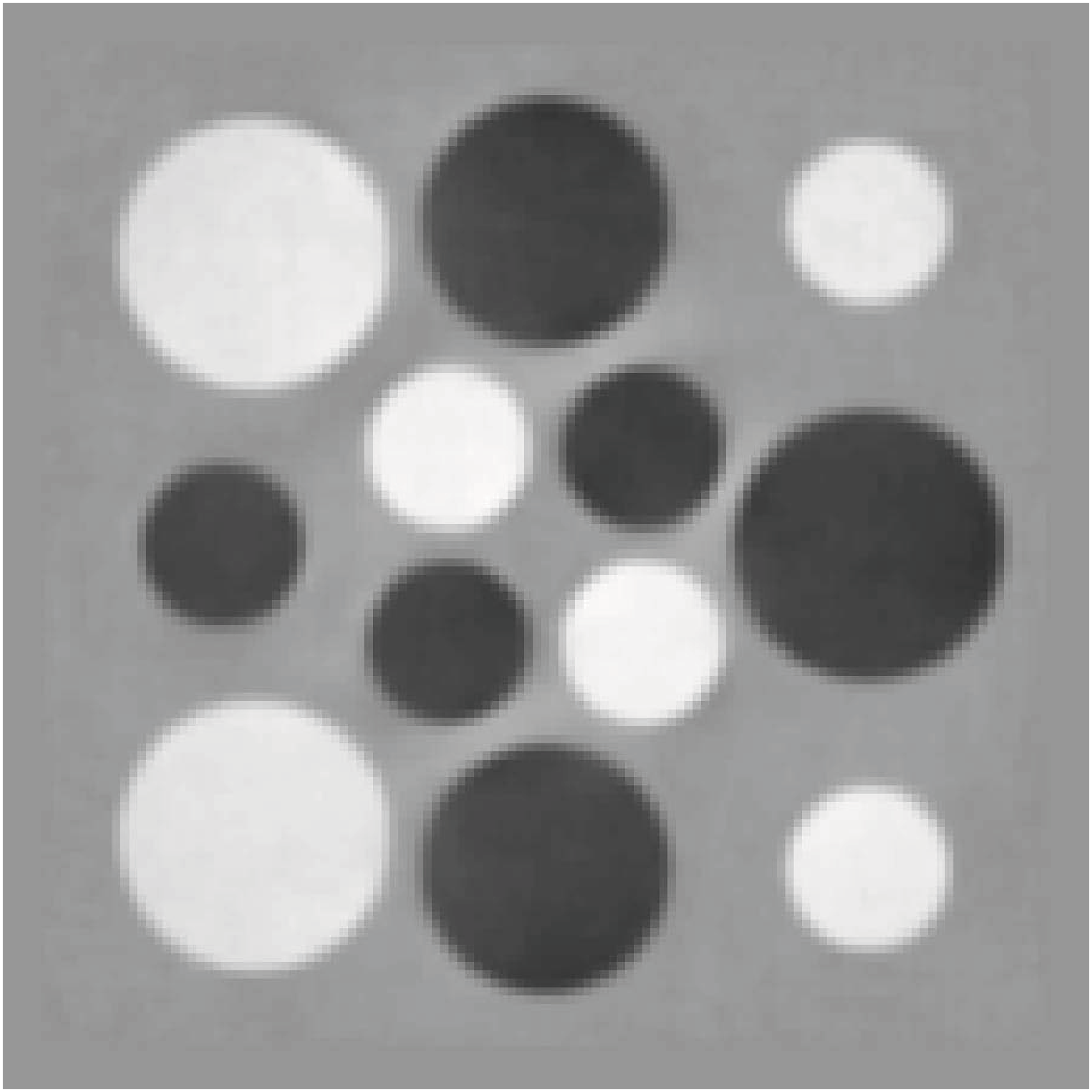} &
\includegraphics[width=1.7in,height=1.7in]{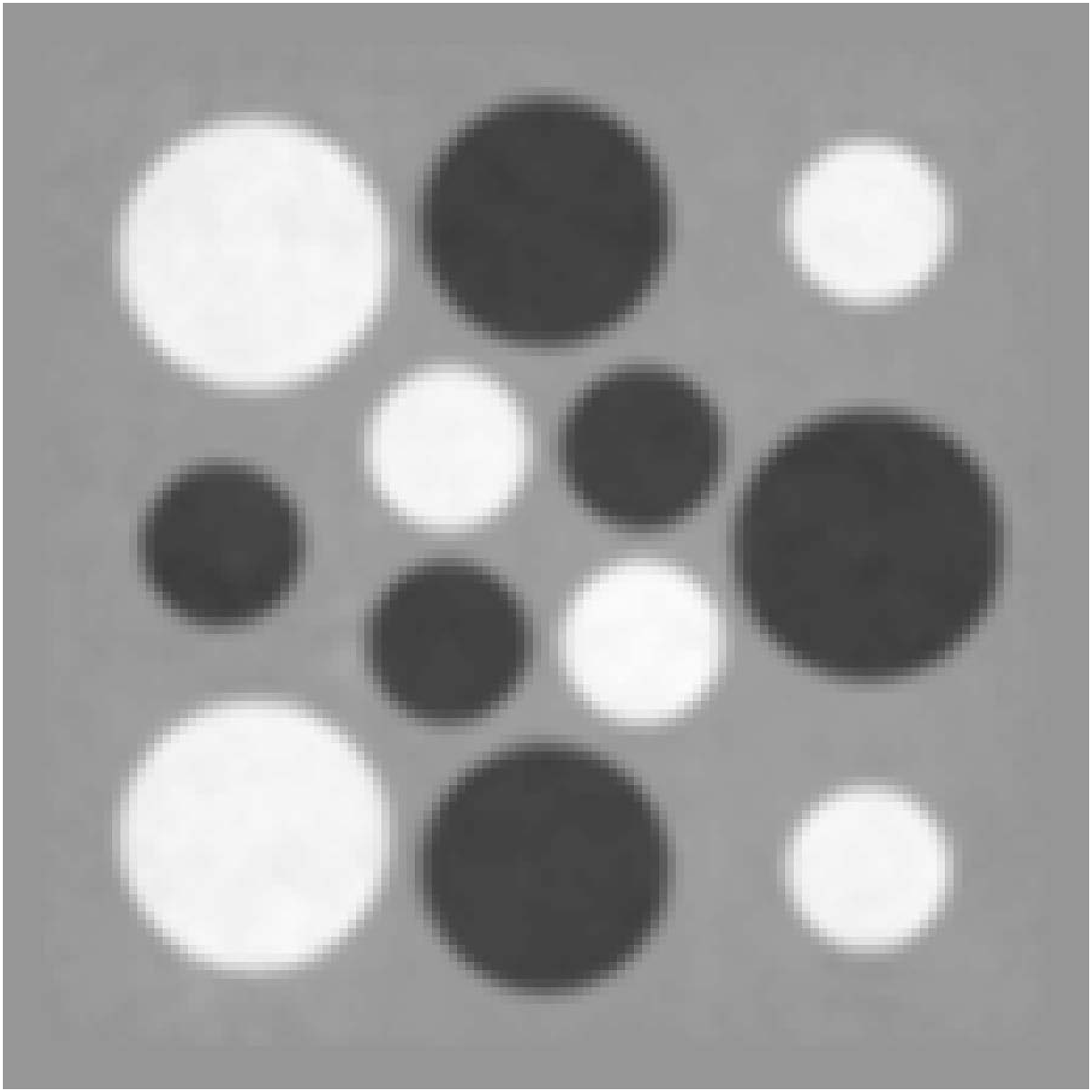}\\
&  & \\
(a) & (b) & (c)
\end{tabular}
\end{center}
\caption{Reconstruction from the data contaminated by a $50\%$ noise
(a)~phantom (b)~iteration~\#0 (c)~iteration~\#1}
\label{F:noisecircgray}
\end{figure}
\begin{figure}[th!]
\begin{center}
\includegraphics[width=2.8in,height=1.1in]{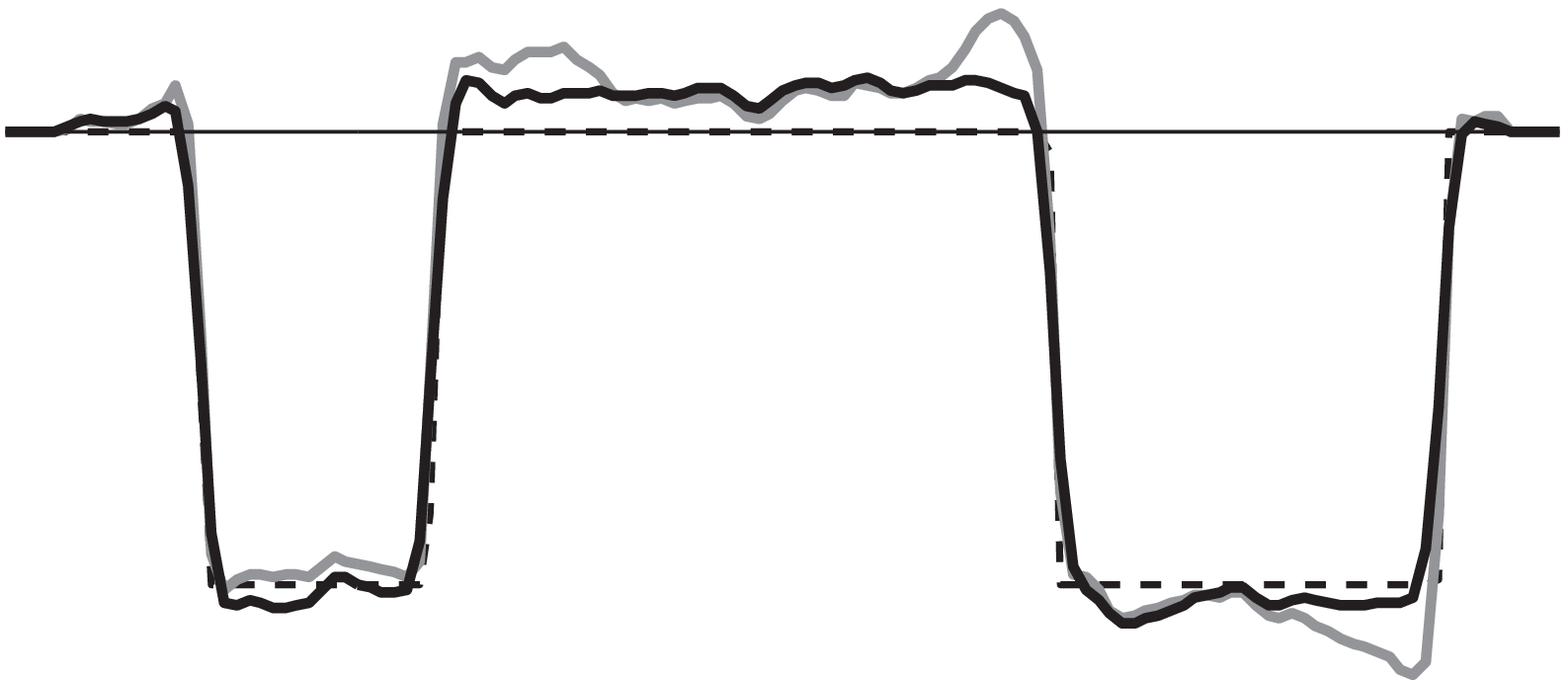}
\end{center}
\caption{Horizontal central cross-section (noisy data): dashed line denotes
the phantom, gray line represents iteration~\#0, thick black solid line
represents iteration~\#1 }
\label{F:noisecircprof}
\end{figure}

 Our phantom (i.e., simulated $\ln\sigma(x)$)
consists of several slightly smoothed characteristic functions of circles,
shown in Figure~\ref{F:accucircgray}(a) and Figure~\ref{F:noisecircgray}(a). (A
more detailed description is presented in the Appendix). Smoothing guarantees
that the phantom is fully resolved on the fine discretization grid we use
during the forward computations, which helps to ensure its high accuracy
(several correct decimal digits). The characteristic functions comprising the
phantom are weighted with weights 1 or -1, so that $\sigma(x)$ varies between
$e$ and $e^{-1}$. Thus, the conductivity deviates far from the initial guess
$\sigma_{0}\equiv1$. Current $I_{1}$ equals $1$ and $-1$ on the right and left
sides of square, respectively; it vanishes on the horizontal sides. Current
$I_{2}$ coincides with $I_{1}$ rotated $90$ degrees counterclockwise.

The simulated sources of the propagating spherical acoustic fronts are centered
on a circle of the diameter slightly larger than the diagonal of the square
domain. There were $256$ simulated transducers uniformly distributed over the
circle. Each transducer produced $257$ spherical fronts of the radii ranging
from $0$ to the diameter of the circle. For each front radius $t_{l}$ and
center $z_{m}$, the perturbed $\sigma$ was modeled, the non-linear direct
$M_{I_{j},I_{k}}(t_{l},z_{m})$, $j,k=1,2$ were computed as explained at the
beginning of this section. In the first of our experiments, these accurate data
were used as a starting point of the reconstruction. In the second experiment,
they were perturbed by a 50\% (in the $L^{2}$ norm) noise.
\begin{figure}[t]
\begin{center}
\begin{tabular} [c]{ccc}
\includegraphics[width=1.7in,height=1.7in]{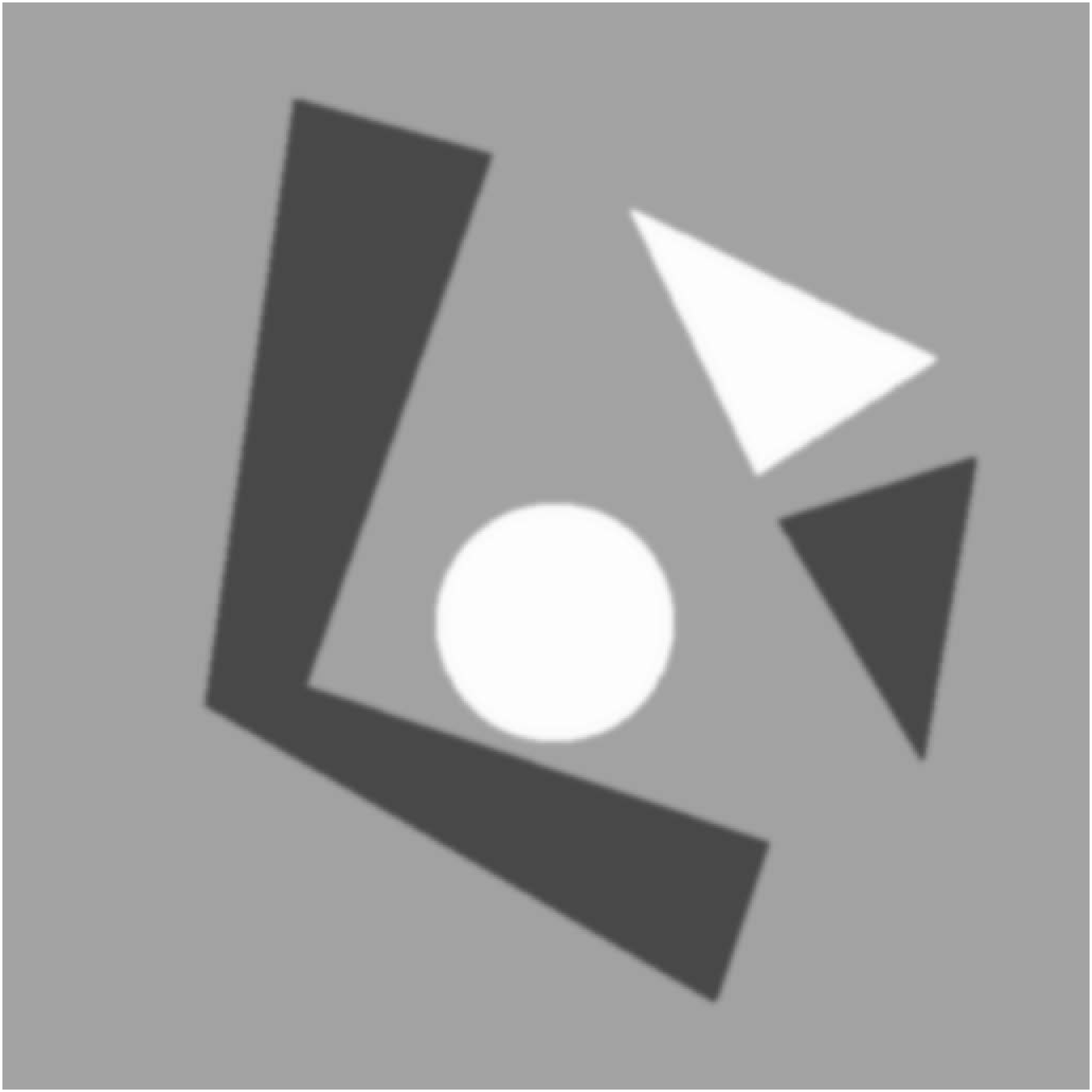} &
\includegraphics[width=1.7in,height=1.7in]{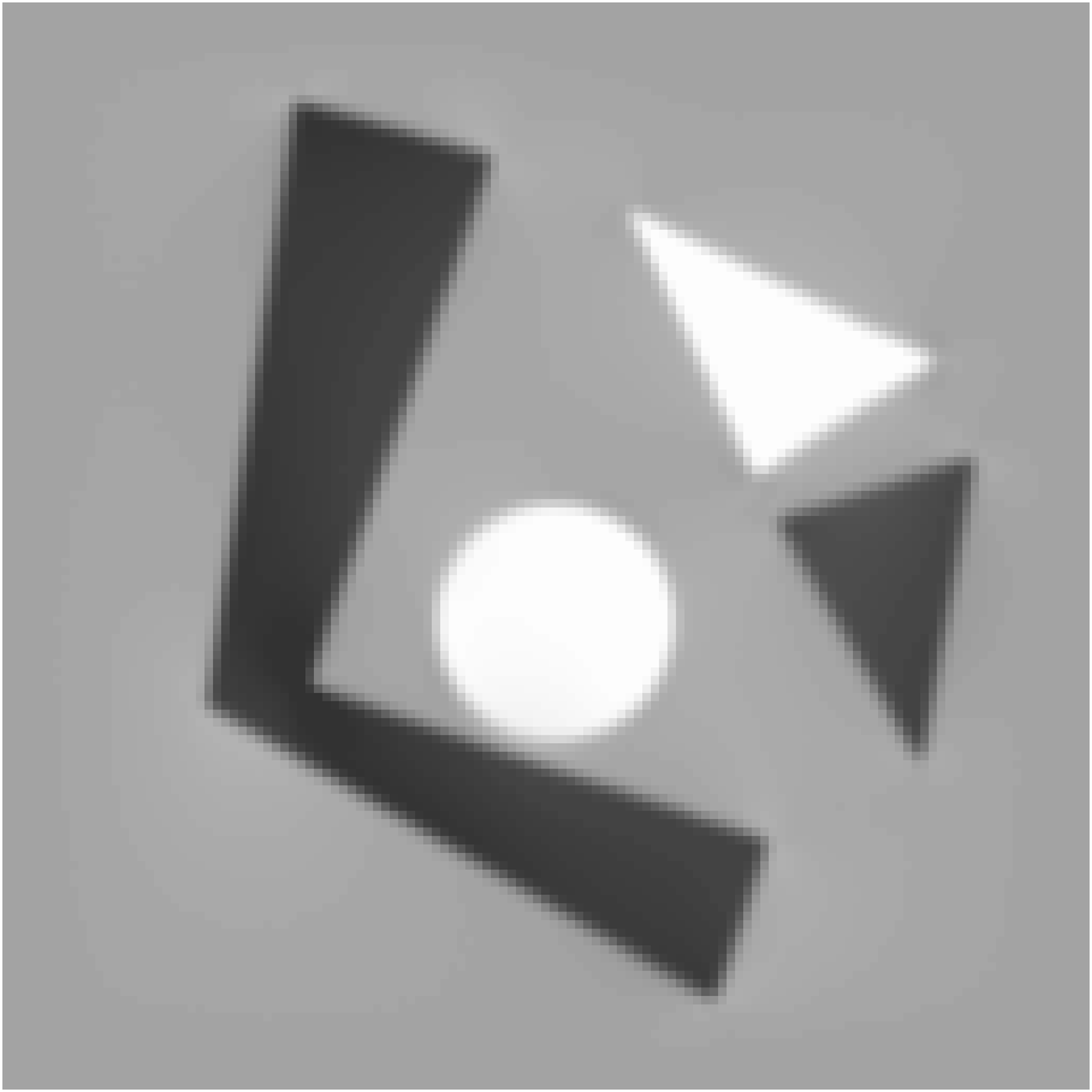} &
\includegraphics[width=1.7in,height=1.7in]{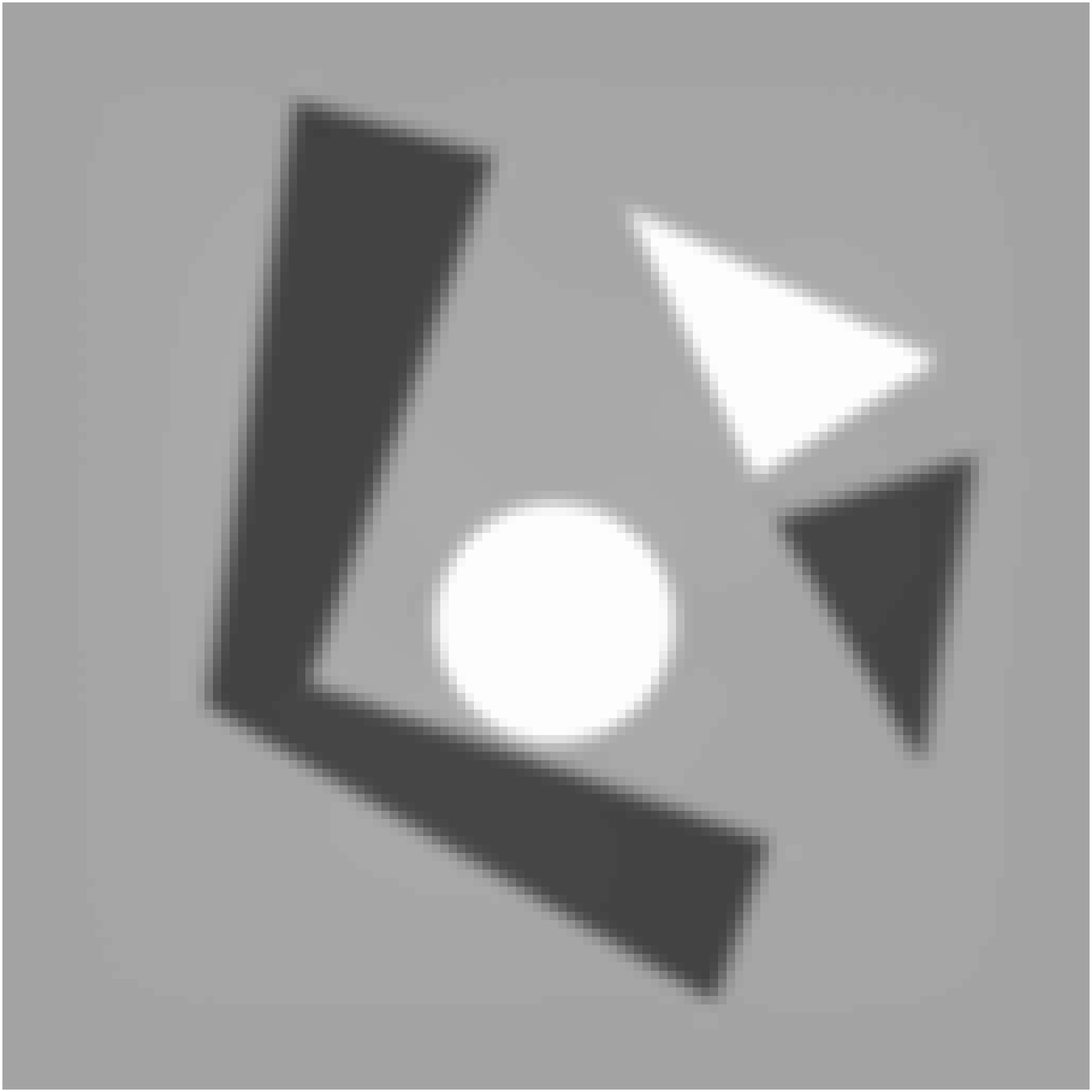}\\
&  & \\ (a) & (b) & (c)
\end{tabular}
\end{center}
\caption{Reconstruction from noiseless data
(a)~phantom (b)~iteration~\#0 (c)~iteration~\#4}
\label{F:accuanglegray}
\end{figure}

The first step of the reconstruction is
synthetic focusing, i.e. finding the values $M_{j,k}(x)$ from $M_{I_{j},I_{k}
}(t,z)$, $j,k=1,2$. In order to give the reader a better feeling of synthetic
focusing, we present in the Figure~\ref{F:focus} a picture of a propagating
spherical acoustic front (part (a)), and an approximation to a delta function
located at the point $(0.2,0.4)$ obtained as a linear combination of such
fronts (part (b)). Figure~\ref{F:focus}(c) shows the same function as in the
part (b) with a modified gray scale that corresponds to the lower $10\%$ of that
function's range, and thus allow one to see small details invisible in part (b).
These figures are provided for demonstration purposes only, since
in our algorithm reconstruction of the values $M_{j,k}(x)$ from $M_{I_{j},I_{k}
}(t,z)$ is done by applying the $2D$ exact filtration
backprojection formula to the latter function (we used the exact reconstruction
formula from \cite{Ku2007}, but other options are also
available).  On a $129\times129$ grid this computation takes a few seconds.
Since the formula is applied to the data containing the derivative of the
delta-function, the differentiation appearing in the TAT inversion formula
(e.g., \cite{AKK,KuKu,FR,Ku2007}) is not needed, and the reconstruction instead
of being slightly unstable, has a smoothing effect (this is why we obtain high
quality images with such high level of noise).

On the second step of the reconstruction,
functions $M_{j,k}^{0}(x)$ are computed using the knowledge of the benchmark
conductivity $\sigma_{0}$, and values of $g_{j,k}(x)$ are obtained by
comparing $M_{j,k}(x)$ and $M_{j,k}^{0}(x)$. Then the first approximation to
$\rho$ (we will call it iteration \#0) is obtained by solving equation
(\ref{E:smartPoisson}). The right hand side of this equation is computed by
finite differences, and then the Poisson equation in a square is solved by the
decomposition in $2D$ Fourier series. The computation is extremely fast due to
the use of the FFT. More importantly, since the differentiation of the data is
followed by the application of the inverse Laplacian, this step is completely
stable (the corresponding pseudodifferential operator is of order zero), and
no noise amplification occurs.
Finally, we attempt to improve the reconstruction by accepting the
reconstructed $\sigma$ as a new benchmark conductivity and
by applying to the data the parametrix
algorithm of the previous section.
We will call this computation iteration \#1.
Figure~\ref{F:accucircgray} demonstrates the result of such reconstruction
from data without noise. Part~(a) of the Figure shows the phantom, parts~(b) and~(c)
present the results of iterations~\#0 and~\#1, on the same gray-level scale.
The profiles of the central horizontal cross-sections of these functions are
shown in Figure~\ref{F:accucircprof}. One can see that even the iteration \#0
produces quite good a reconstruction; iteration~\#1 removes some of the
artifacts, and improves the shape of circular inclusions. For the convenience
of the reader we summarize the parameters of this simulation in the Appendix.

Figures~\ref{F:mij}, \ref{F:noisecircgray} and \ref{F:noisecircprof} present
the results of the reconstruction from noisy data. In this simulation we used
the phantom from the previous example, and we added to the data $50\%$ (in
$L^2$ norm) noise. The first step of the reconstruction (synthetic focusing) is
illustrated by Figure~\ref{F:mij}. Parts~(a) and~(c) of this Figure show
accurate values of the functionals $M_{1,1}(x)$ and $M_{1,2}(x)$. Parts~(b)
and~(d) present the reconstructed values of these functionals obtained by
synthetic focusing. One can see the effect of smoothing mentioned earlier in
this section: the level of  noise in the reconstructions is much lower than the
level of noise in the simulated measurements. The images reconstructed from
$M_{i,j}(x)$  on the second step are presented in Figures~\ref{F:noisecircgray}
and \ref{F:noisecircprof}. The meaning of the images is the same as of those in
Figures~\ref{F:accucircgray} and \ref{F:accucircprof}. The level of noise in
these images is comparable to that in the reconstructed  $M_{i,j}$'s. To
summarize, our method can reconstruct high quality images from the data
contaminated by a strong noise since the first step of the method is an
application of a smoothing operator, and the second step uses the parametrix.

Finally, Figure~\ref{F:accuanglegray} shows
reconstruction of a phantom containing objects with corners. The phantom is
shown in the part~(a) of the figure, part~(b) demonstrates iteration~\#0, and
part~(c) presents the result of the iterative use of the parametrix method
described in the previous section (iteration~\#4 is shown).

\section{Reconstruction in $3D$}\label{S:3drec}

Let us now consider the reconstruction problem in $3D$. The $3D$ case is very
important from the practical point of view, since propagation of electrical
currents is essentially three-dimensional. Indeed, unlike X-rays or
high-frequency ultrasound, currents cannot be focused to stay in a
two-dimensional slice of the body.
However, while successful $3D$ reconstructions were reported~\cite{Cap},
the theoretical foundations of the $3D$ case have not been
completed yet, due to some analytic difficulties
arising in other approaches.
In contrast, the present approach easily
generalizes to $3D$, and leads to a fast, efficient, and robust reconstruction
algorithm.

We will assume that three different currents $I_{j},j=1,2,3$ are used, and
that the boundary values of the corresponding potentials $w_{j},$ $j=1,2,3$
are measured on $\partial\Omega.$ Similarly to the $2D$ case presented in
Section \ref{S:formulation}, by perturbing the medium with a perfectly focused
acoustic beam (no matter whether such measurements are real or synthesized)
one can recover at each point $x$ within $\Omega$ the values of the
functionals $M_{i,j}(x),$ $i,j=1,2,3,$ where, as before,
\begin{equation}
M_{i,j}(x)=\sigma(x)\nabla w_{i}(x)\cdot\nabla w_{j}(x). \label{E:M3d}
\notag \end{equation}

Our goal is to reconstruct conductivity $\sigma(x)$ from $M_{i,j}(x).$ As
before, we will assume that $\sigma(x)$ is a perturbation of a known benchmark
conductivity $\sigma_{0}(x),$ i.e. $\sigma(x)=\sigma_{0}(x)(1+\varepsilon
\rho(x)),$ and that the values of  potentials $w_{j}(x)$ are the perturbations
of known potentials $u_{j}(x)$ corresponding to $\sigma_{0}(x):$
\begin{equation}
w_{j}(x)=u_{j}(x)+\varepsilon v_{j}(x)+o(\varepsilon).
\notag \end{equation}
Now functionals $M_{j,k}(x)$ are related to the known unperturbed values
$M_{j,k}^{0}(x)$ and measured perturbations $g_{j,k}(x)$ by equations
(\ref{E:M}) and (\ref{E:M0}).

As it was done in Section \ref{S:reconstruction}, we introduce vector fields
$U_{j}=\sqrt{\sigma_{0}}\nabla u_{j}$ and $W_{j}=\sqrt{\sigma}\nabla
(u_{j}+\varepsilon v_{j})=U_{j}+\varepsilon V_{j},$ and proceed to derive the
following six equations:
\begin{align*}
U_{1}\cdot V_{1} &  =g_{1,1}/2\\
U_{2}\cdot V_{2} &  =g_{2,2}/2\\
U_{3}\cdot V_{3} &  =g_{3,3}/2\\
U_{1}\cdot V_{2}+U_{2}\cdot V_{1} &  =g_{1,2}\\
U_{1}\cdot V_{3}+U_{3}\cdot V_{1} &  =g_{1,3}\\
U_{2}\cdot V_{3}+U_{3}\cdot V_{2} &  =g_{2,3}.
\end{align*}

One can obtain a useful approximation to $\rho(x)$ by assuming $\sigma_{0}=1$,
and by selecting unperturbed currents so that the potentials $u_{j}(x)=x_{j}$.
Then, by repeating derivations of Section \ref{S:constant} one obtains the
following three formulas
\begin{equation}
\left\{
\begin{array}
[c]{c}
\left(  \frac{\partial^{2}}{\partial x_{1}^{2}}+\frac{\partial^{2}}{\partial
x_{2}^{2}}\right)  \rho=\frac{1}{2}\left(  \frac{\partial^{2}}{\partial
x_{1}^{2}}-\frac{\partial^{2}}{\partial x_{2}^{2}}\right)  (g_{2,2}
-g_{1,1})-2\frac{\partial^{2}}{\partial x_{1}\partial x_{2}}g_{1,2}\\
\left(  \frac{\partial^{2}}{\partial x_{1}^{2}}+\frac{\partial^{2}}{\partial
x_{3}^{2}}\right)  \rho=\frac{1}{2}\left(  \frac{\partial^{2}}{\partial
x_{1}^{2}}-\frac{\partial^{2}}{\partial x_{3}^{2}}\right)  (g_{3,3}
-g_{1,1})-2\frac{\partial^{2}}{\partial x_{1}\partial x_{3}}g_{1,3}\\
\left(  \frac{\partial^{2}}{\partial x_{2}^{2}}+\frac{\partial^{2}}{\partial
x_{3}^{2}}\right)  \rho=\frac{1}{2}\left(  \frac{\partial^{2}}{\partial
x_{2}^{2}}-\frac{\partial^{2}}{\partial x_{3}^{2}}\right)  (g_{3,3}
-g_{2,2})-2\frac{\partial^{2}}{\partial x_{2}\partial x_{3}}g_{2,3}
\end{array}
\right.  \label{E:3dsystem}
\end{equation}
We notice that by using the first of the above equations one can compute an
approximation to $\rho(x)$ by solving a set of $2D$ Poisson equations (one for
each fixed value of $x_{3})$, since boundary values of $\rho(x)$ are equal to
0. This leads to a slice-by-slice $3D$ reconstruction, which is based only on
values of $g_{1,1},$ $g_{2,2}$ and $g_{1,2}$, and therefore can be done by
using a single pair of currents.

One can get better images by using all three currents and doing a fully $3D$
reconstruction. Namely, summing the equations (\ref{E:3dsystem}) yields the
values of $2\Delta\rho$ in the left hand side. Then one can solve the $3D$
Poisson equation with the zero boundary conditions to recover the
conductivity.

One can expect that, as in $2D$, this approach would work well for $\sigma(x)$
close to $\sigma_{0}=1$. However, as demonstrated by our numerical experiments
presented in Section \ref{S:3d}, the results remain quite accurate when
$\sigma(x)$ varies significantly across $\Omega$. Moreover, a simple fixed
point iteration based on the repeated use of formulas (\ref{E:3dsystem})
exhibits a rapid convergence to the correct image.

\section{Numerical examples in $3D$} \label{S:3d}

\begin{figure}[t]
\begin{center}
\begin{tabular}
[c]{ccc}
\includegraphics[width=1.7in,height=1.7in]{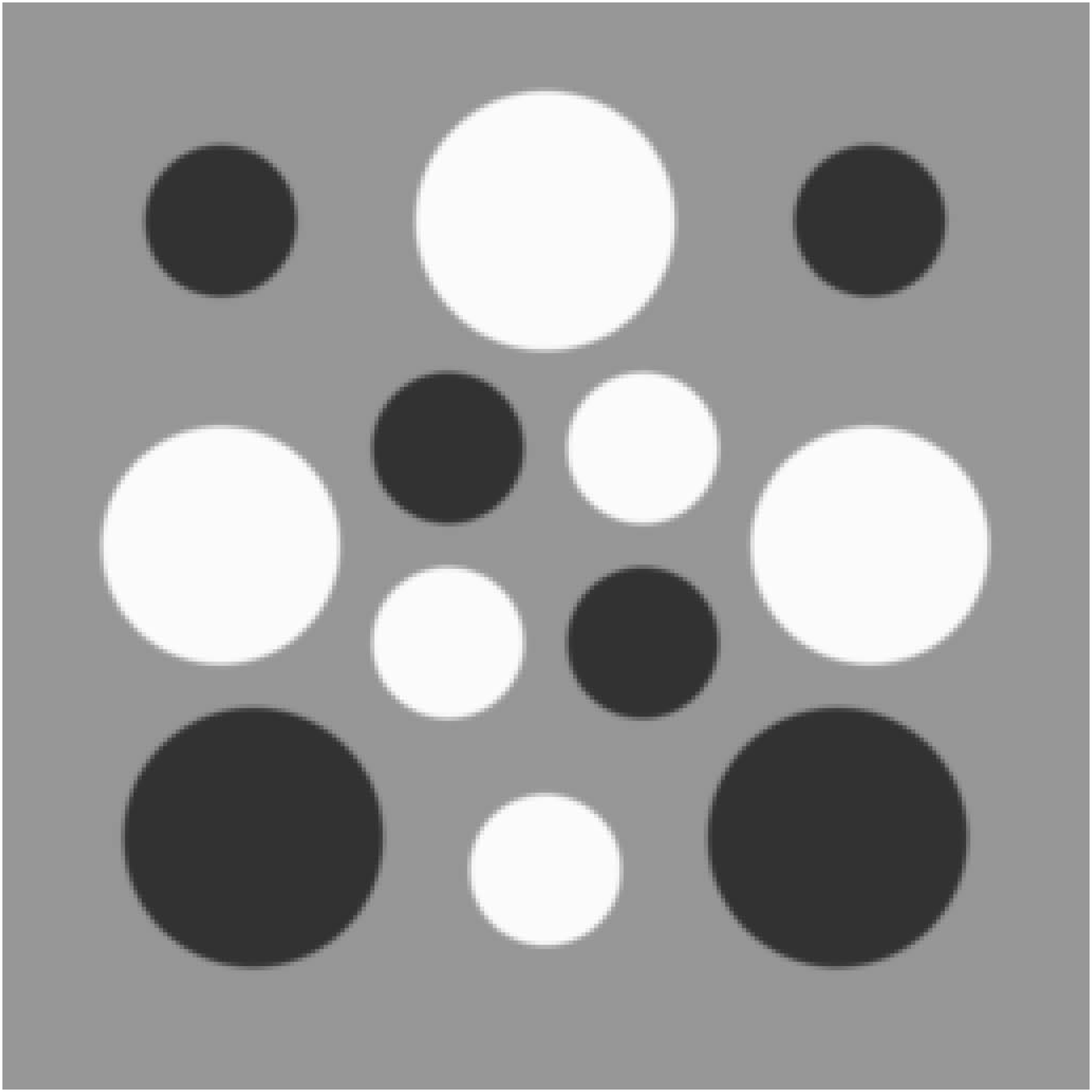} &
\includegraphics[width=1.7in,height=1.7in]{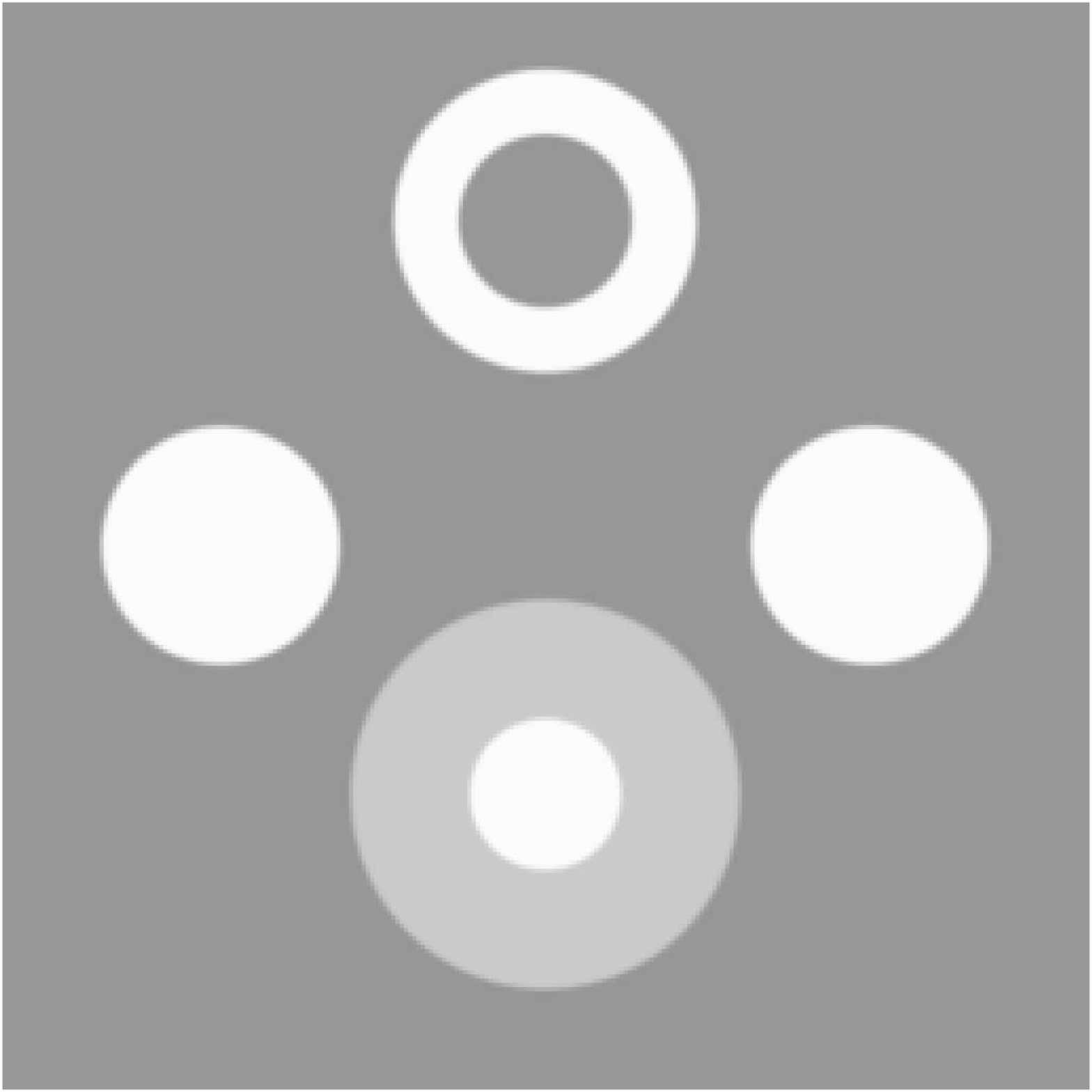} &
\includegraphics[width=1.7in,height=1.7in]{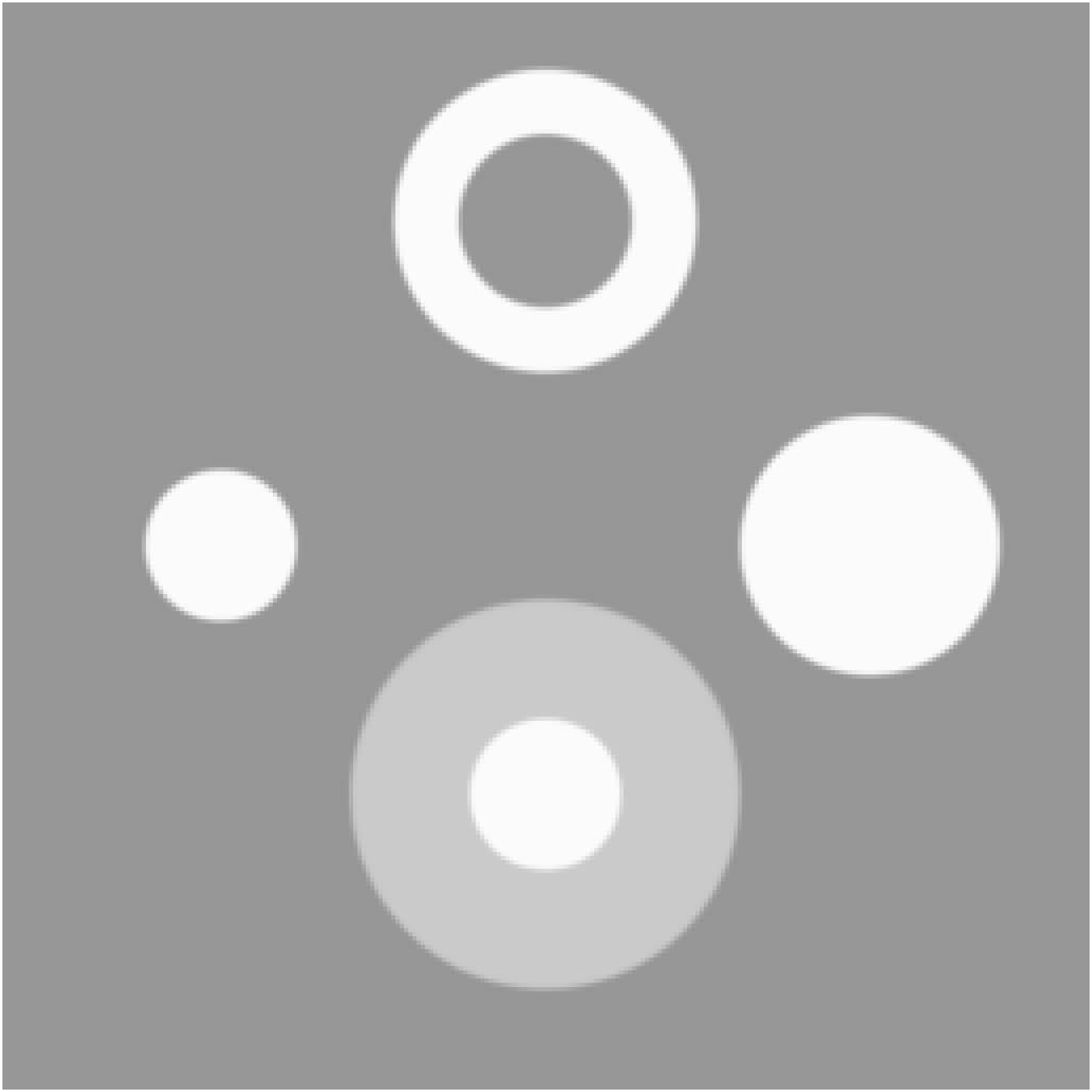}\\
\includegraphics[width=1.7in,height=1.7in]{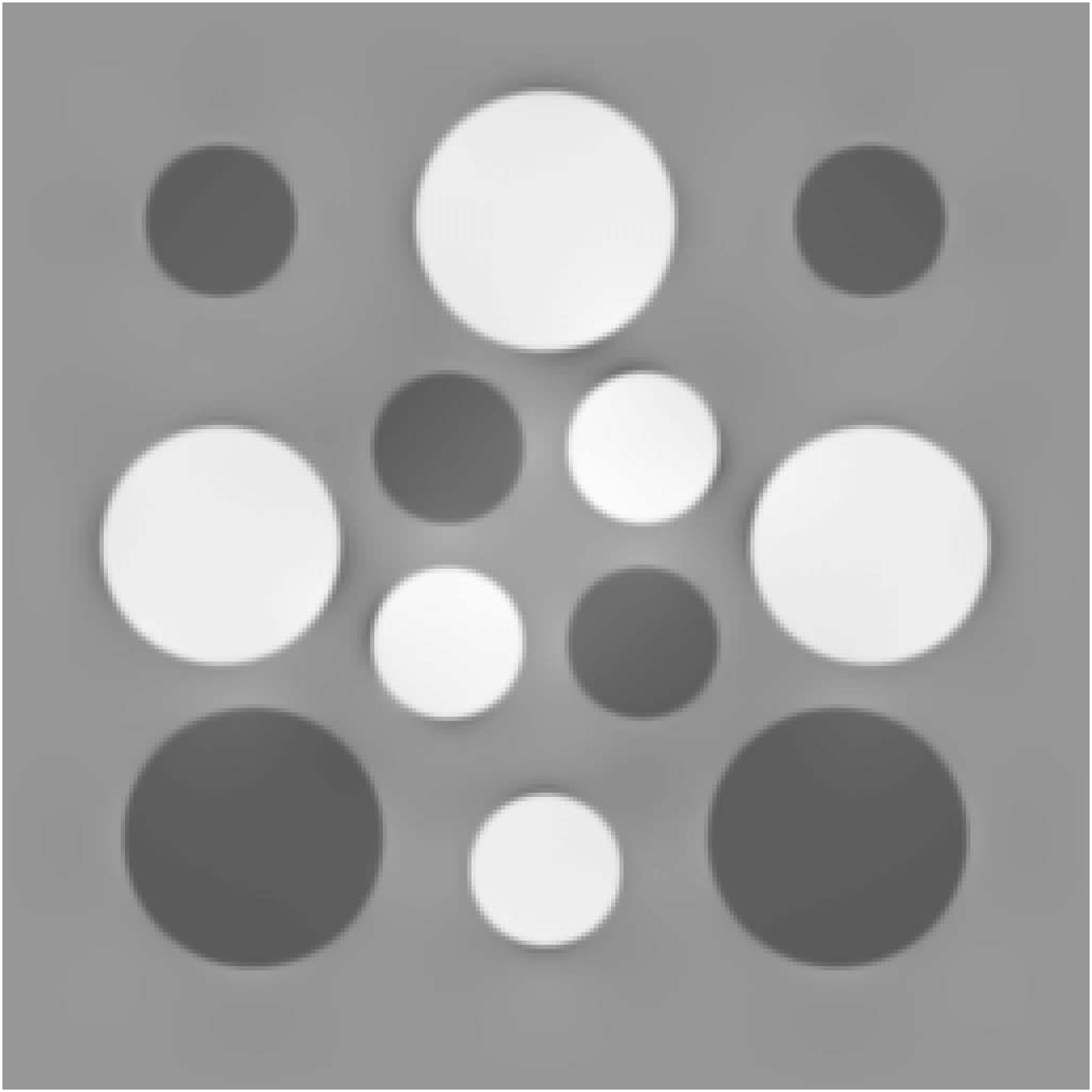} &
\includegraphics[width=1.7in,height=1.7in]{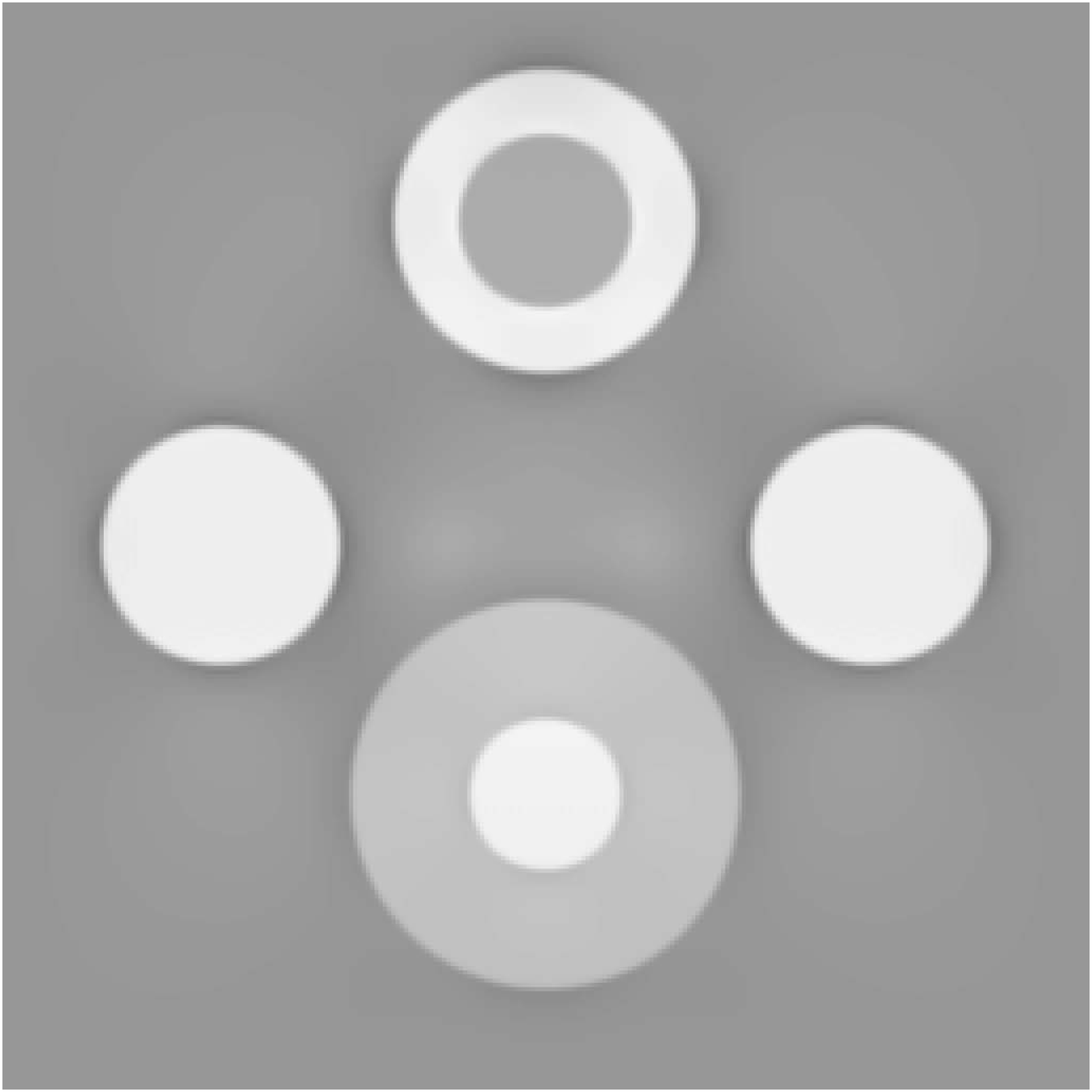} &
\includegraphics[width=1.7in,height=1.7in]{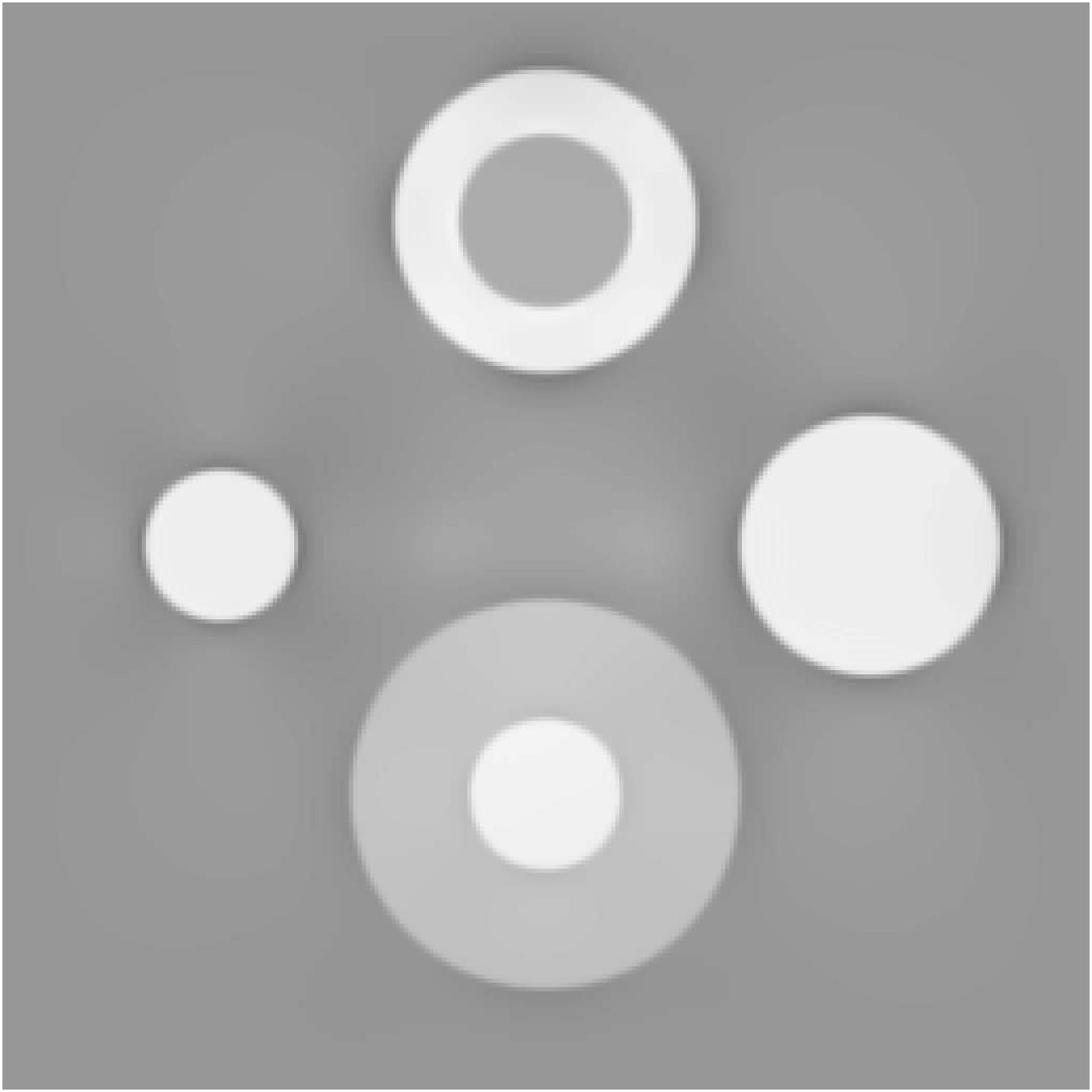}\\
\includegraphics[width=1.7in,height=1.7in]{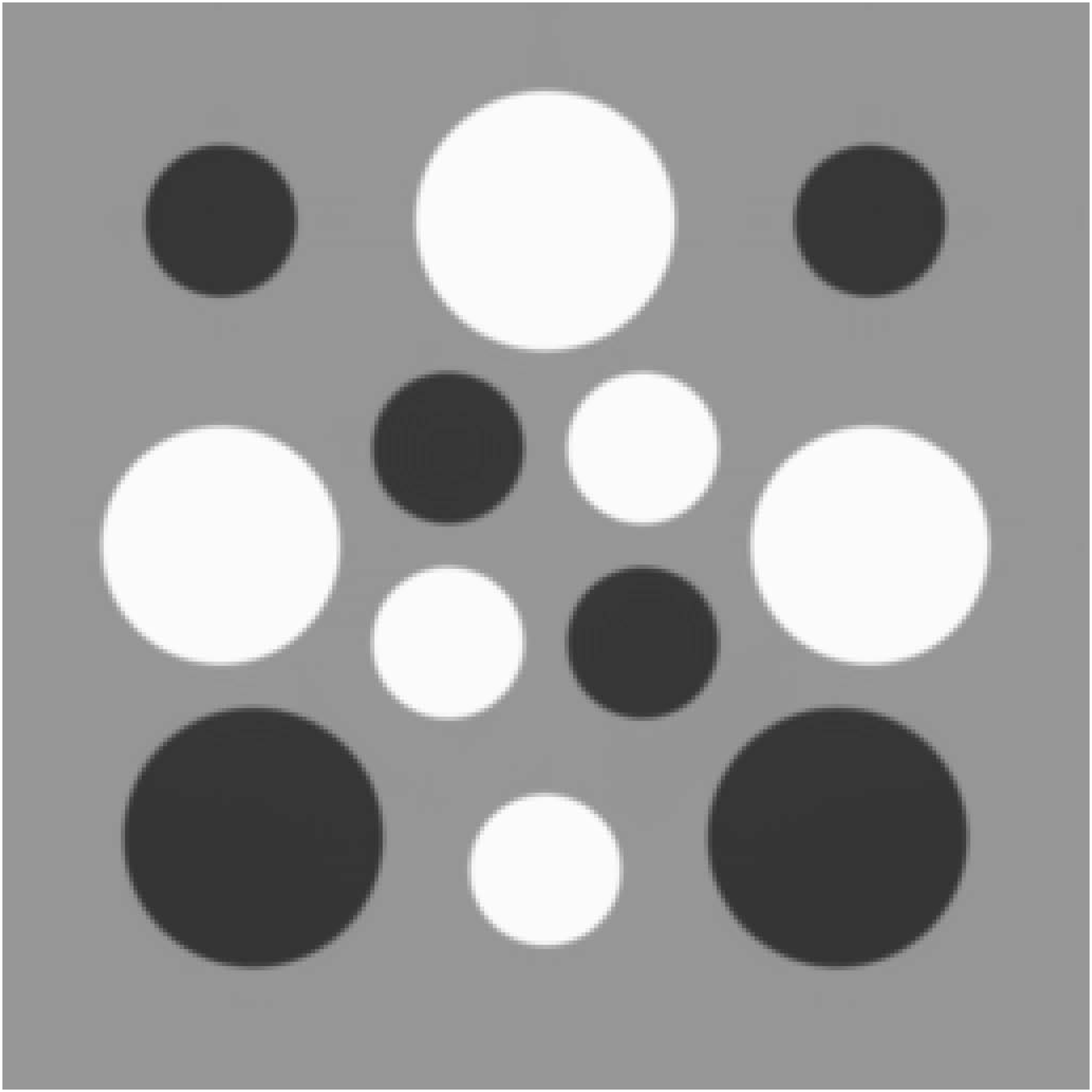} &
\includegraphics[width=1.7in,height=1.7in]{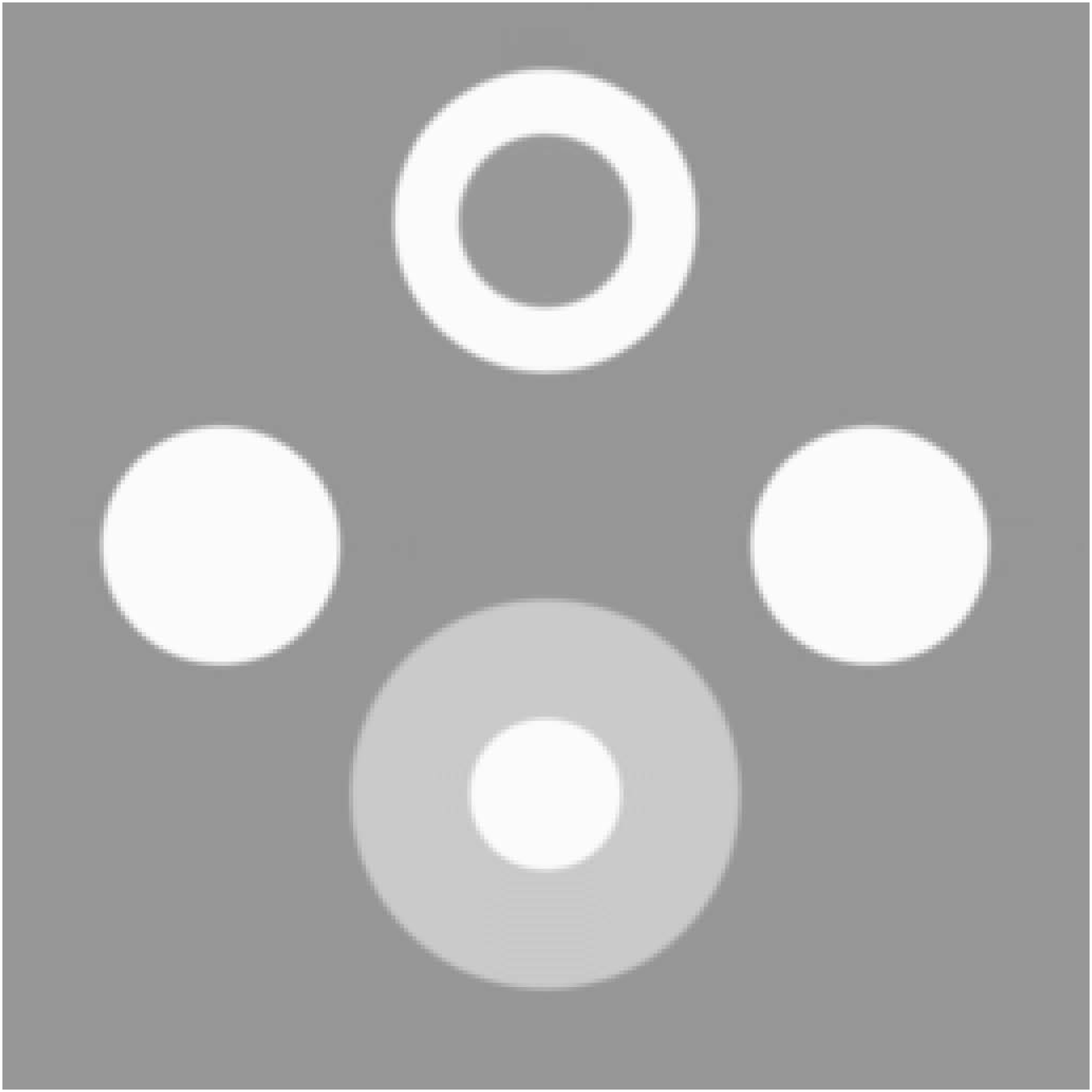} &
\includegraphics[width=1.7in,height=1.7in]{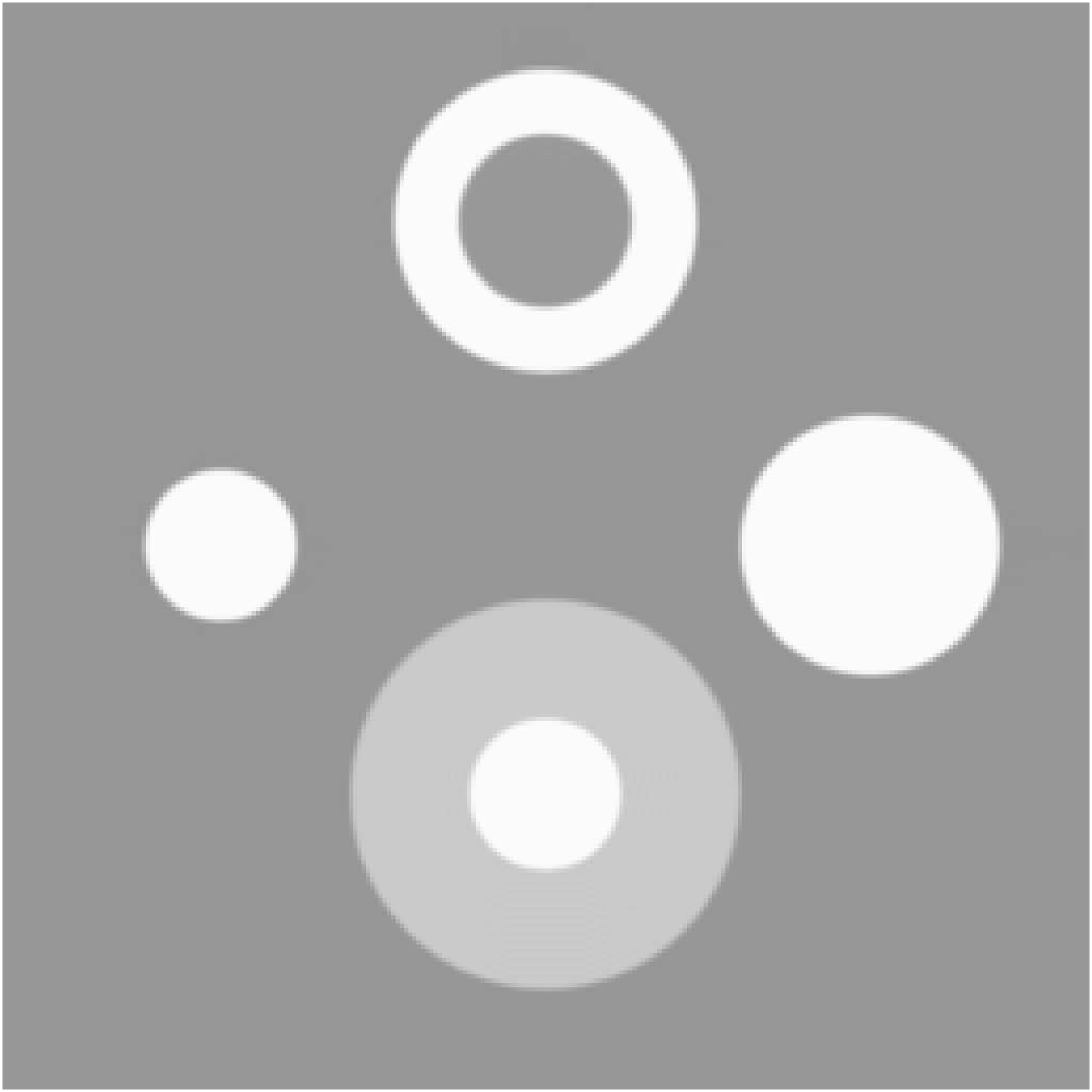}\\
(a) & (b) & (c)
\end{tabular}
\end{center}
\caption{$3D$ Reconstruction from noiseless data. First row: phantom
(a)~$Ox_1x_2$~cross section
(b)~$Ox_1x_3$~cross section
(c)~$Ox_2x_3$~cross section.
Second row: iteration~\#0; Third row: iteration~\#4}
\label{F:accu3Dbw}
\end{figure}

\begin{figure}[th!]
\begin{center}
\includegraphics[width=2.8in,height=1.1in]{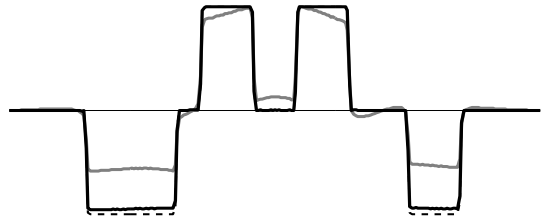}
\end{center}
\caption{Diagonal cross-section (noiseless data): dashed line denotes the
phantom, gray line represents iteration~\# 0, thick black solid line represents
iteration~\#4 }
\label{F:accuprof}
\end{figure}

In this section we present results of $3D$ reconstructions from simulated data.
Unfortunately, a complete modeling of the forward problem in $3D$ (i.e.
computation of the perturbations corresponding to the
propagating acoustic spherical fronts) would require solution of
$\mathcal{O}(n^3)$ $3D$ divergence equations. This task is computationally
too expensive. Therefore, unlike in our $2D$ simulations, we resort to modeling
the values of the functionals $M_{i,j}(x)$ on a $257 \times 257 \times 257$
Cartesian grid, using formulas~(\ref{E:M3d}).
These values correspond to the data that would be measured if perfectly
focused, infinitely small perturbations were applied to the conductivity.
Thus, in this section we only test the second step of our reconstruction
techniques. However, as mentioned before, if the real data were available,
the first step (synthetic focusing) could be done by applying any of the
several available stable versions of thermoacoustic inversion, and the
feasibility of this step was clearly demonstrated in the $2D$ sections of this
paper, as well as in~\cite{KuKuAET}.

In our first simulation we used noiseless values of $M_{i,j}(x)$ and
reconstructed the conductivity on a $257 \times 257 \times 257$ grid.
The first row of Figure~\ref{F:accu3Dbw} shows three $2D$ cross-sections of a
$3D$~phantom. The result of approximate inversion (using three currents,
as described in Section~\ref{S:3drec}) is presented in the second row of the
figure. Finally, the last row shows the result of iterative use of formulas
(\ref{E:3dsystem}), where $\rho$ now represents the difference between the previous
and the updated approximations to the conductivity. The third row demonstrates
iteration~\#4. In addition, Figure~\ref{F:accuprof} shows the trace along a
diagonal cross section in $Ox_1x_2$ plane (that corresponds to the diagonals of
images presented in the column~(a) of Figure~\ref{F:accu3Dbw}). We summarize
the details of this simulation in the Appendix.

\begin{figure}[t!]
\begin{center}
\begin{tabular}
[c]{ccc}
\includegraphics[width=1.7in,height=1.7in]{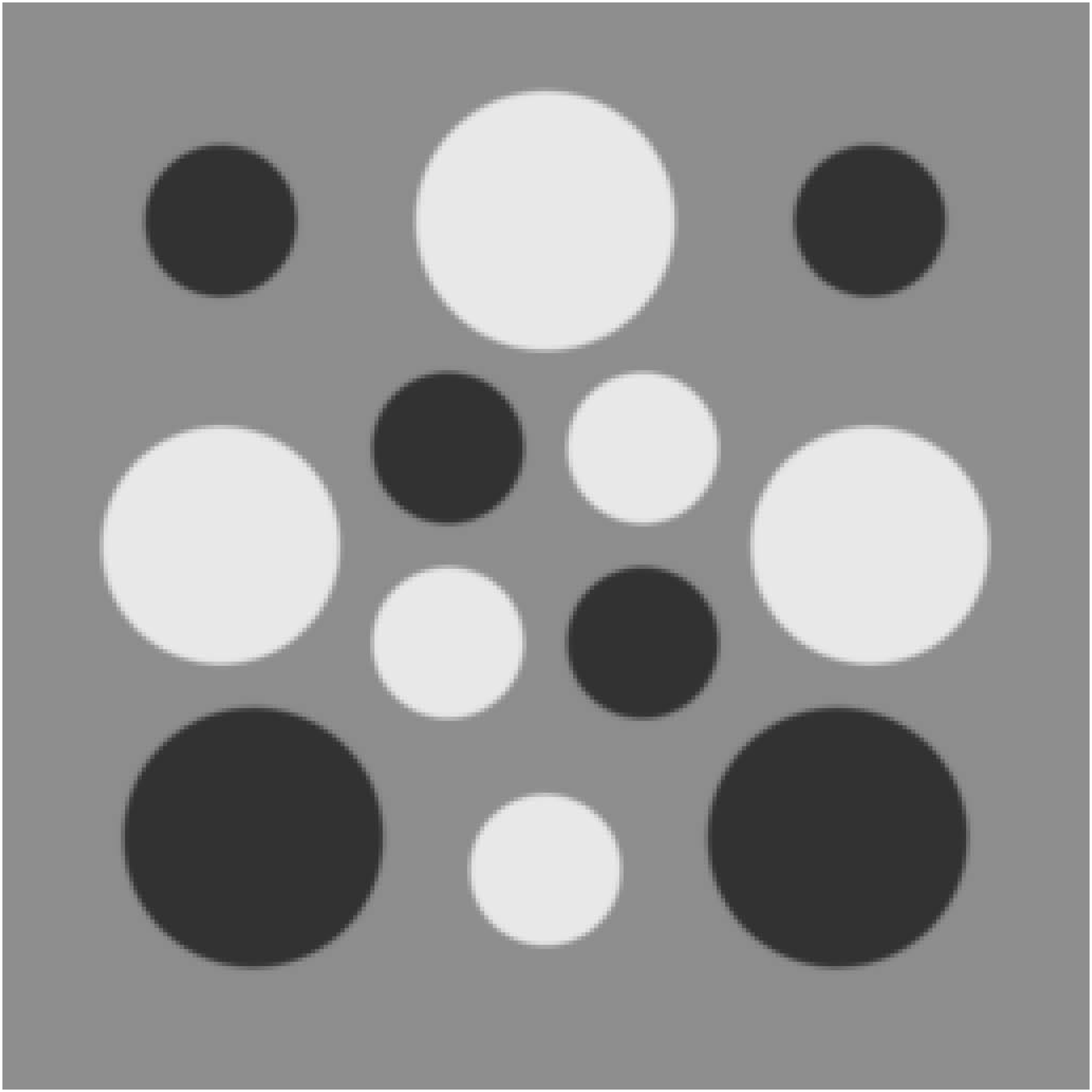} &
\includegraphics[width=1.7in,height=1.7in]{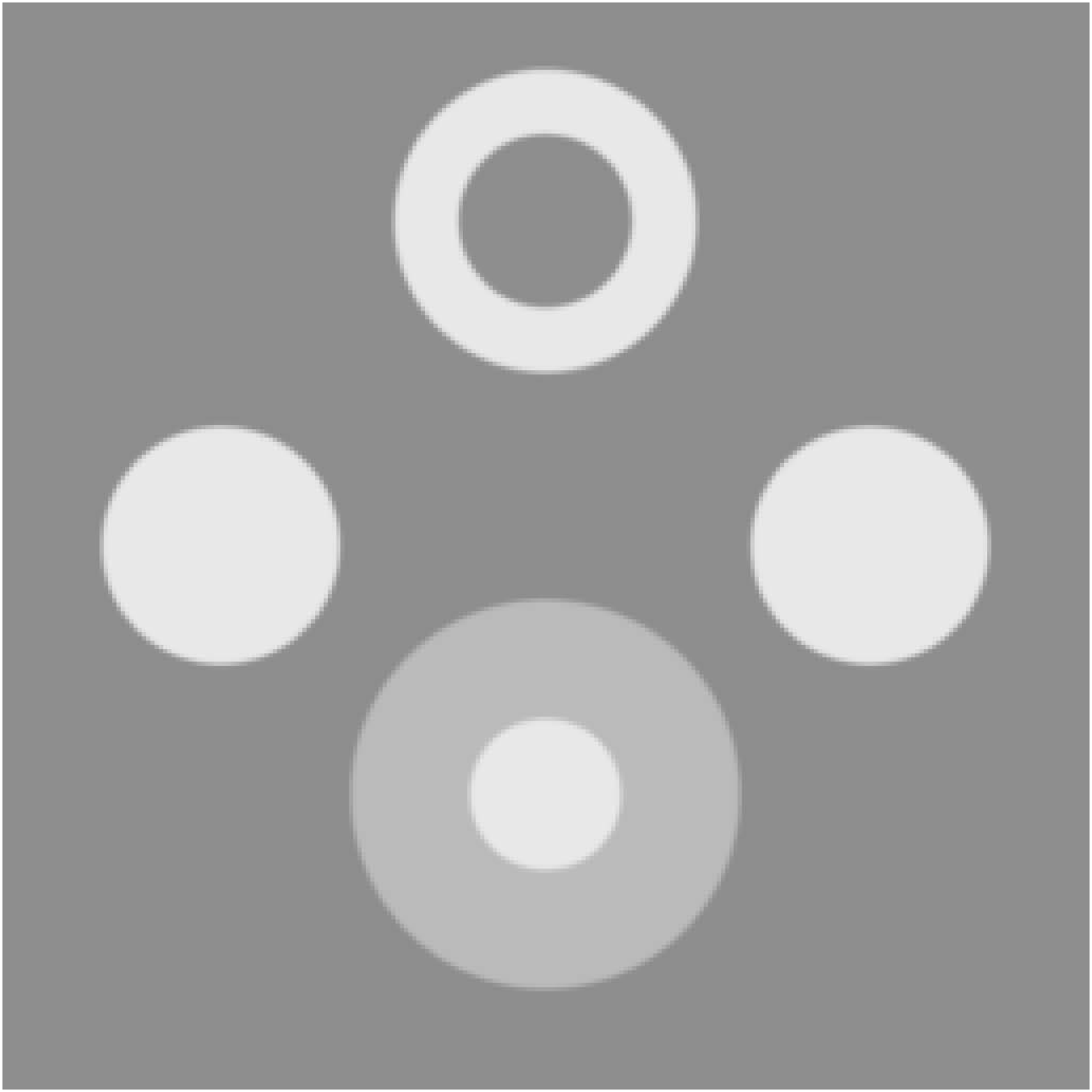} &
\includegraphics[width=1.7in,height=1.7in]{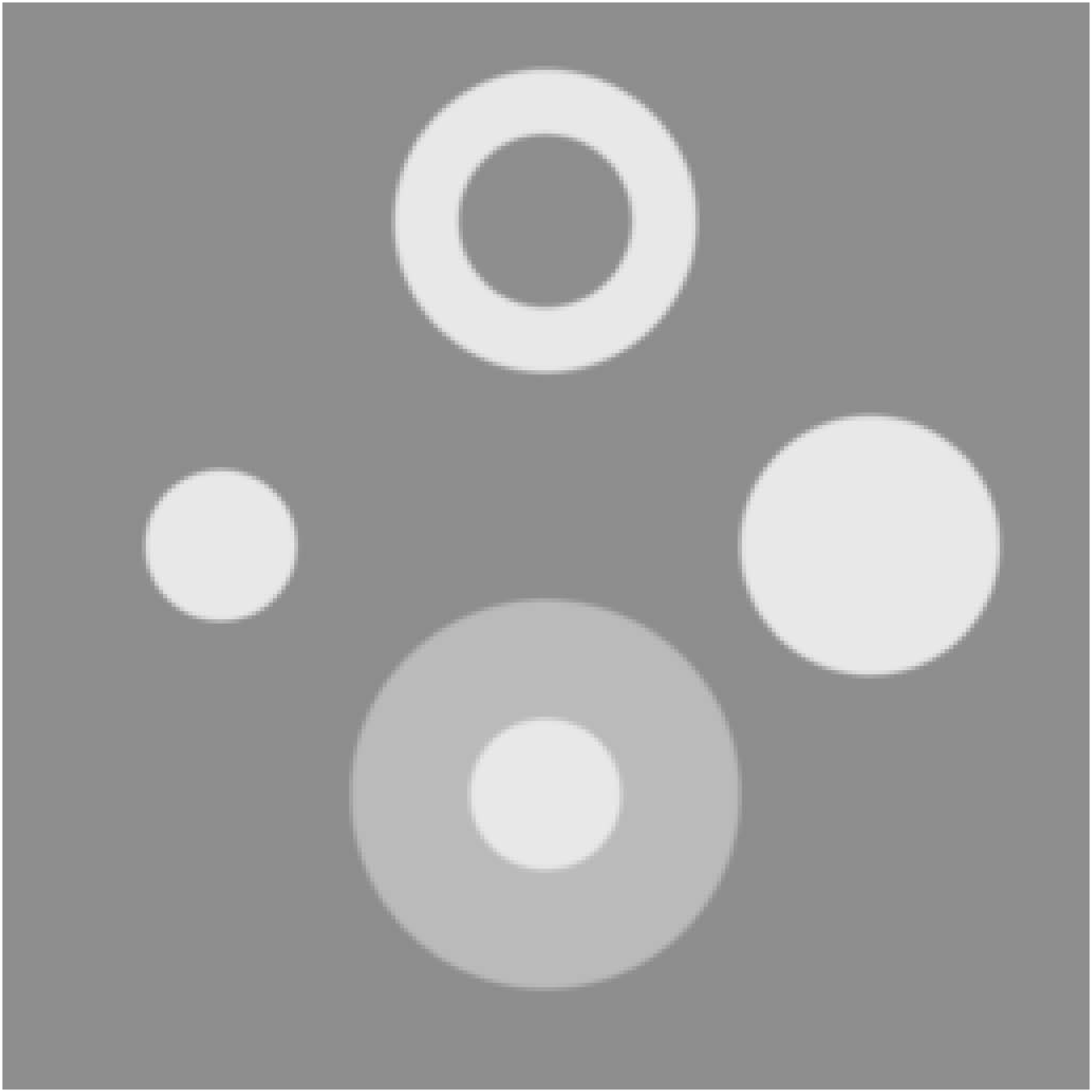}\\
\includegraphics[width=1.7in,height=1.7in]{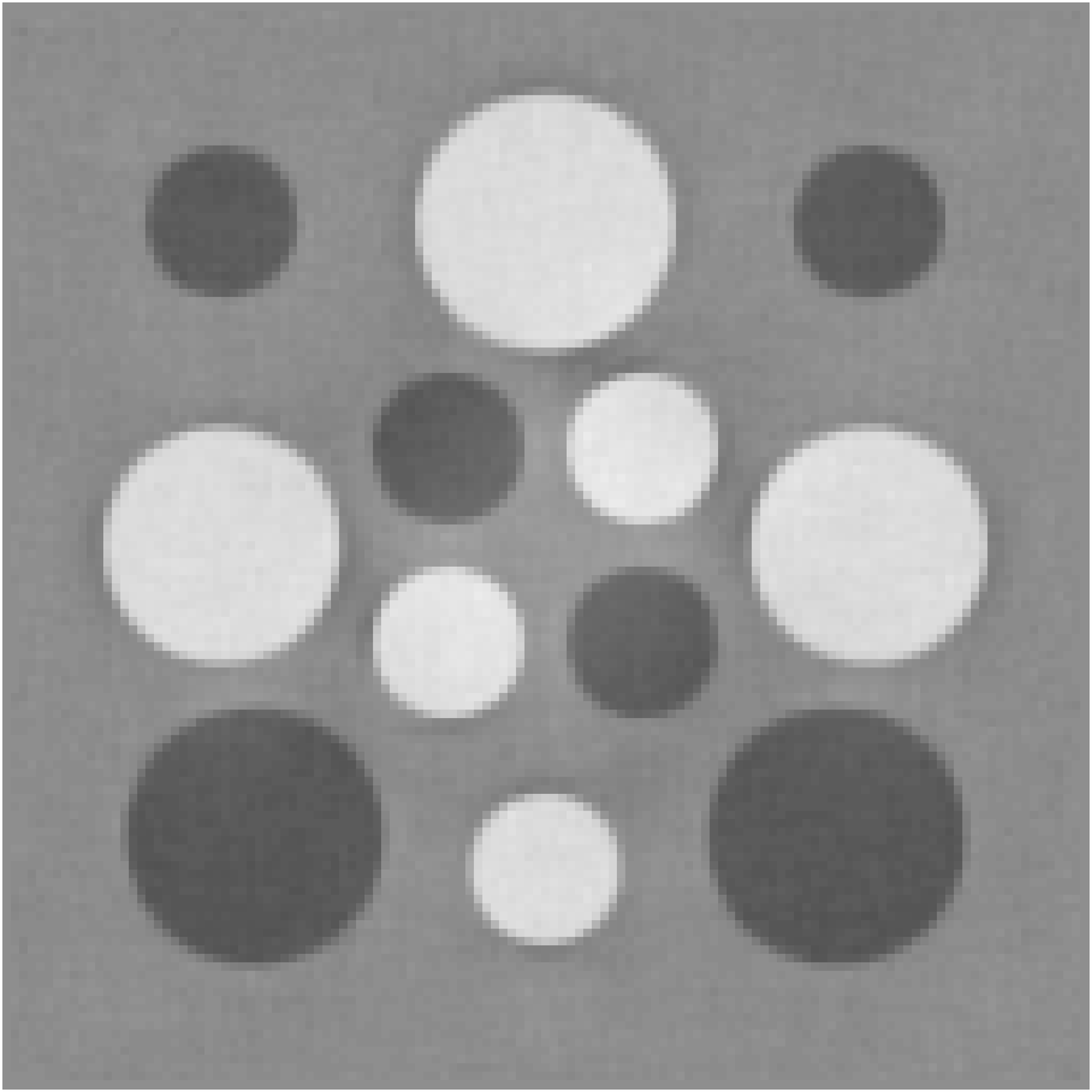} &
\includegraphics[width=1.7in,height=1.7in]{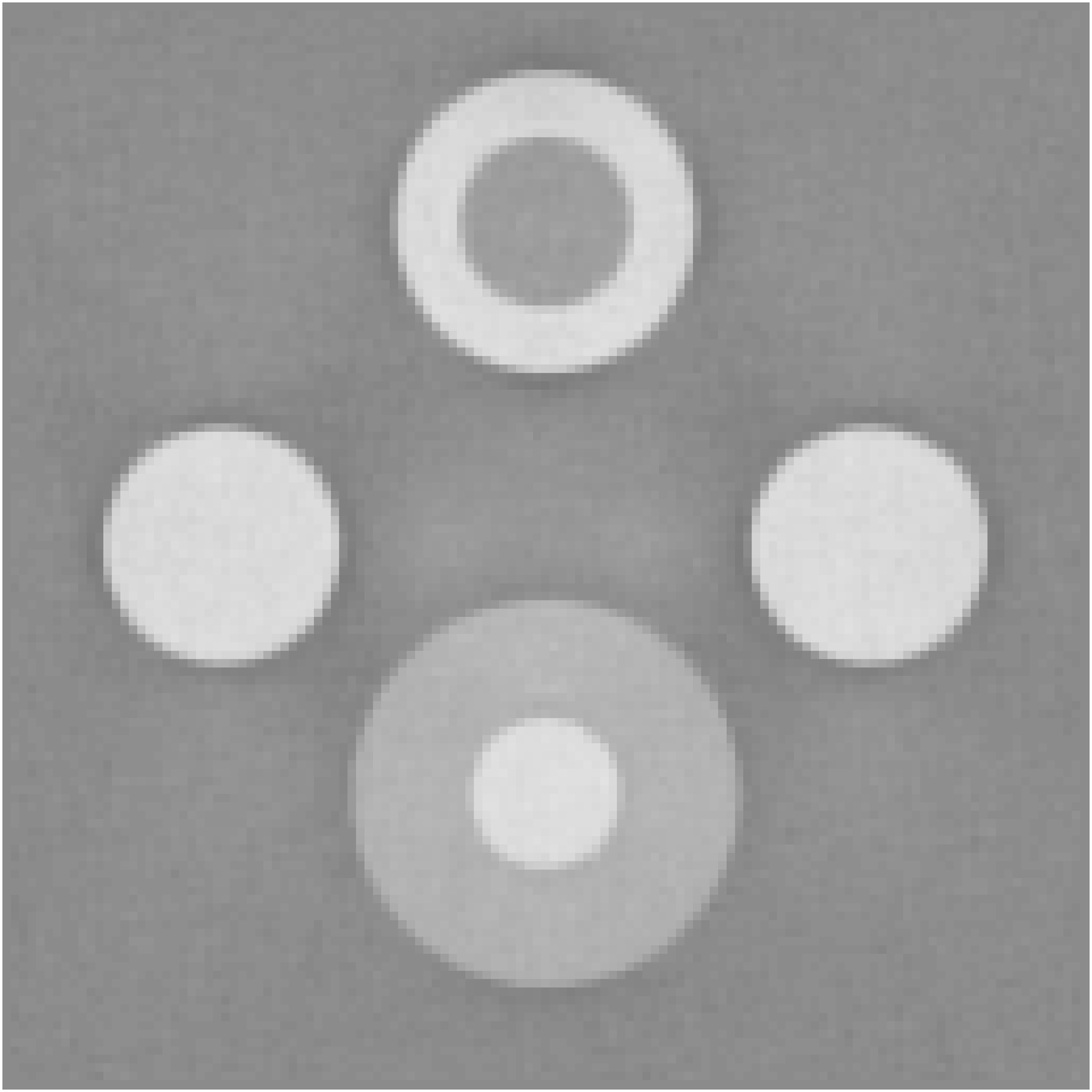} &
\includegraphics[width=1.7in,height=1.7in]{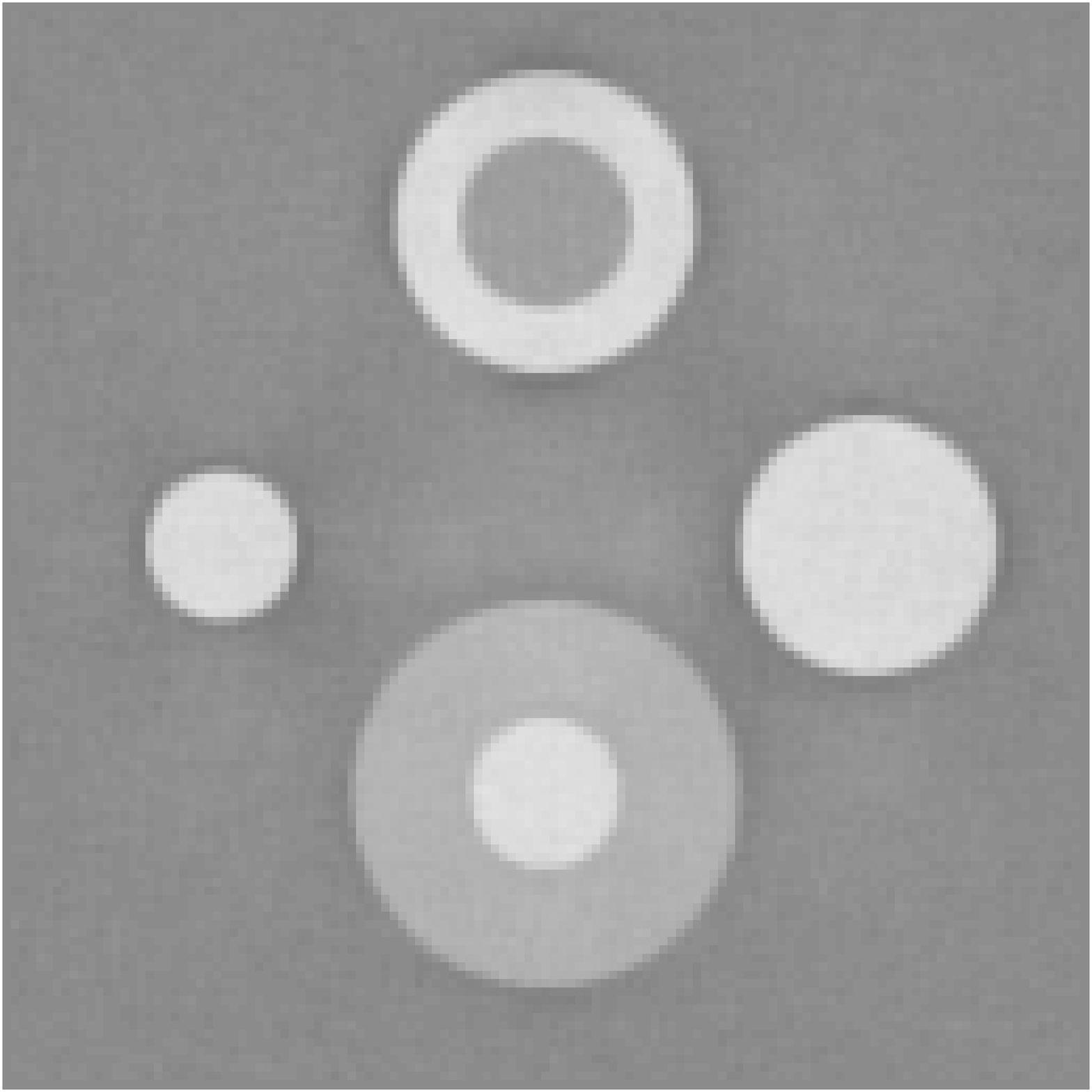}\\
\includegraphics[width=1.7in,height=1.7in]{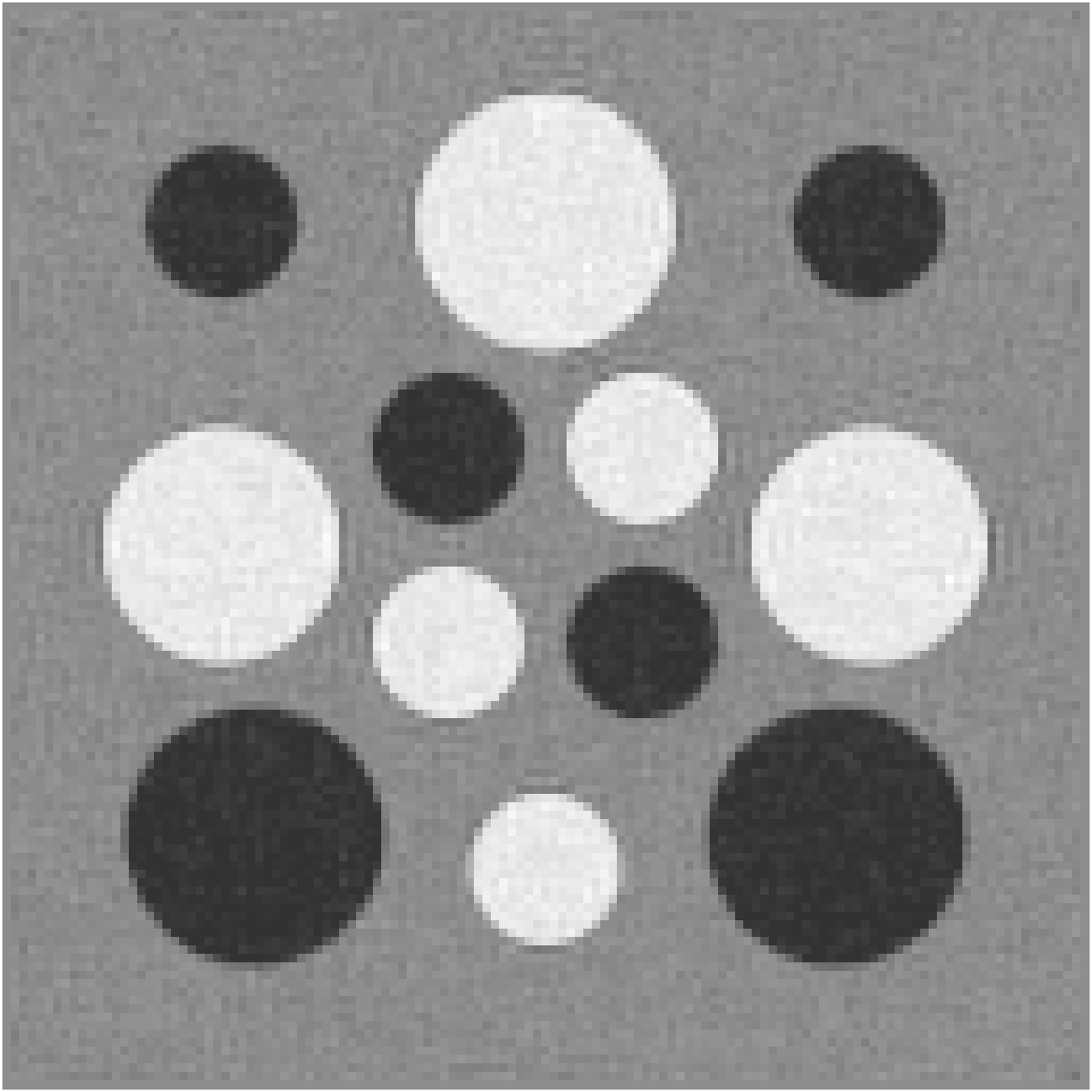} &
\includegraphics[width=1.7in,height=1.7in]{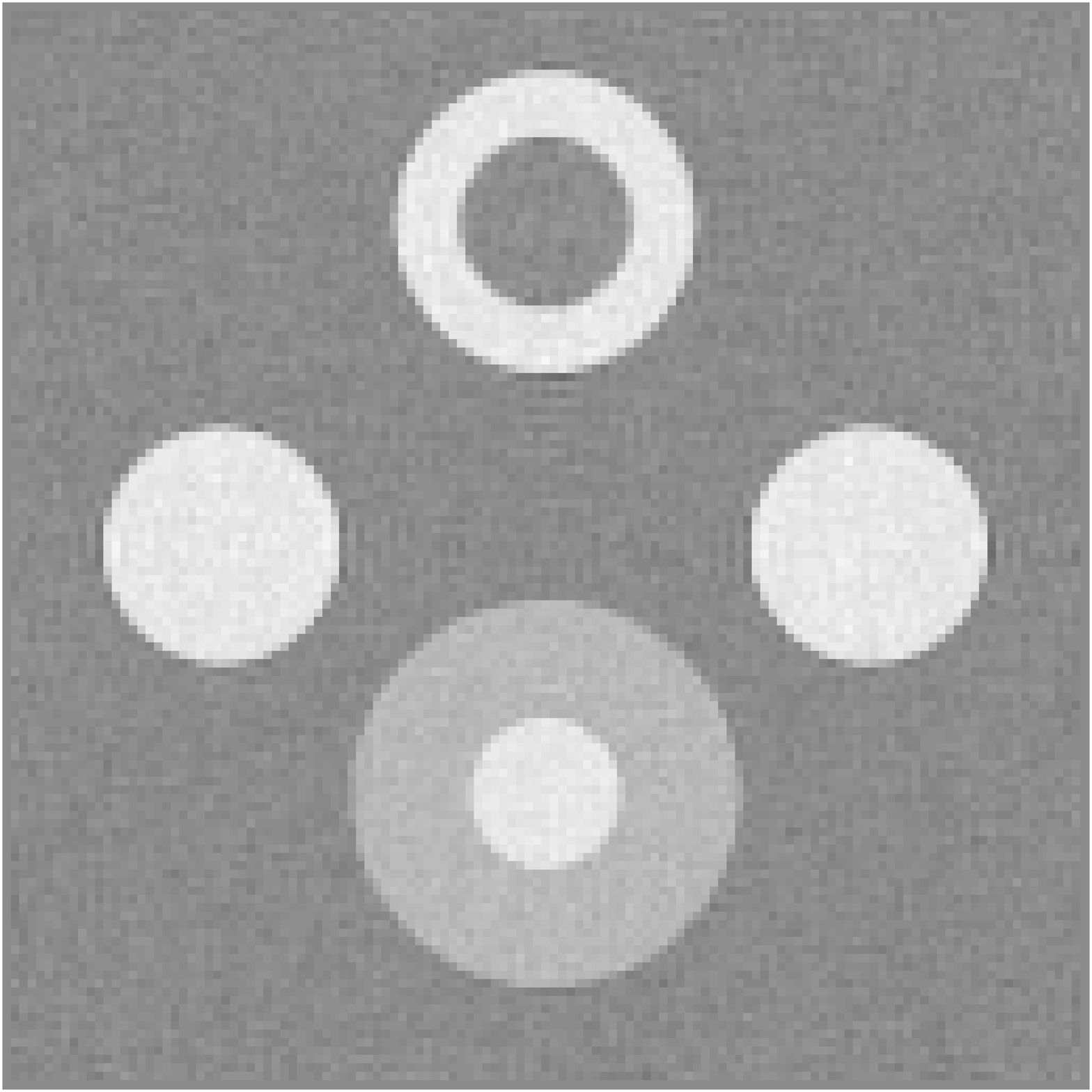} &
\includegraphics[width=1.7in,height=1.7in]{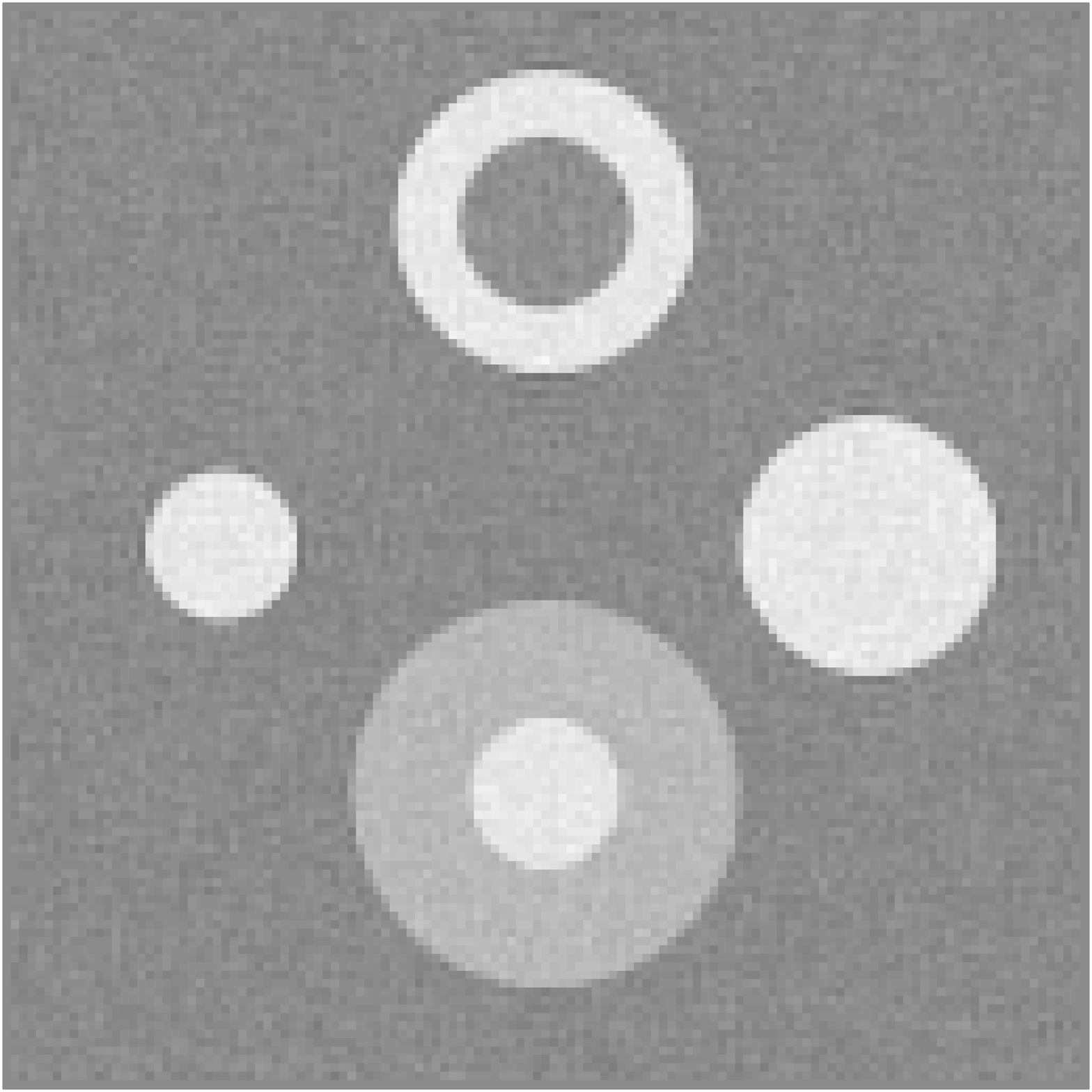}\\
(a) & (b) & (c)
\end{tabular}
\end{center}
\caption{$3D$ Reconstruction from noisy data on a coarser grid. First row: phantom
(a)~$Ox_1x_2$~cross section (b)~$Ox_1x_3$~cross section (c)~$Ox_2x_3$ cross section.
Second row: iteration~\#0; Third row: iteration~\#4}
\label{F:noisegray}
\end{figure}

\begin{figure}[h!]
\begin{center}
\includegraphics[width=2.8in,height=1.2in]{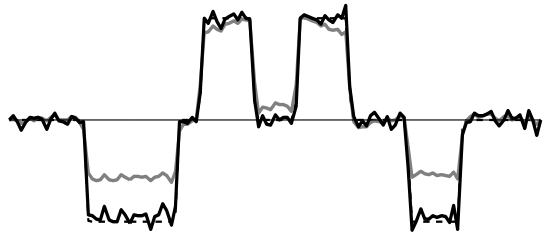}
\end{center}
\caption{Diagonal cross-section of reconstructions obtained from noisy data
on a coarser gird: dashed line denotes
the phantom, gray line represents iteration~\# 0, thick black solid line
represents iteration~\#4 }
\label{F:noiseprof}
\end{figure}

In our second $3D$ experiment we utilized the same phantom, but as a data
used only a subset of the values of $M_{i,j}$ corresponding to a coarser
$129 \times 129 \times 129$ grid; the latter coarse grid was also used to
discretize the reconstructed conductivity. We also added to the data
a $10\%$ (in $L^2$ norm) noise. Figure~\ref{F:noisegray} presents the
cross-sections of a $3D$ phantom and the reconstructions obtained using three
currents, on the same gray-level scale. The meaning of the subfigures
is the same as of those in Figure~\ref{F:accu3Dbw}. Finally,
Figure~\ref{F:noiseprof} shows the trace along the diagonal cross sections
of the images in the $Ox_1x_2y$ plane.

In both these examples iteration~\#0 yields good qualitative reconstruction
of the conductivity in spite the fact that the latter varies from $e^{-1}$ to
$e^1$, and thus differs strongly from the benchmark guess $\sigma_0 = 1$.
The subsequent iterations demonstrate fast convergence to the
correct values of $\sigma(x)$.

\section{Final remarks and conclusions} \label{S:remarks}

We have shown that the proposed algorithm works stably and yields quality
reconstructions of the internal conductivity. It does not require physical
focusing of ultrasound waves and replaces it with the synthetic focusing
procedure, which can be implemented using one of the known thermoacoustic
imaging inversion methods (e.g., time reversal or inversion formulas). Under
appropriate smoothness conditions on the conductivity, our analysis leads to
the proof of local uniqueness and stability of the reconstruction. However,
since this conclusion has been already made in $2D$ in \cite{Cap,bon}, we only
presented a sketch of the proof.

Some additional remarks:
\begin{enumerate}
\item Using the propagating spherical fronts of the type considered in this text
(equation~(\ref{E:pulse})) is advantageous since in this case the synthetic focusing
is a smoothing operator, and thus the whole reconstruction procedure is more stable with
respect to errors than the one that starts with focused data.

\item Reconstructions can be done with a single, two, or (in $3D$) three
currents. A single current procedure was the one we used initially in $2D$
\cite{KuKuAET,Oberw}. It works, but requires solving a transport equation for
the conductivity. When such a procedure is used, errors arising due to the
noise and/or underresolved interfaces tend to propagate along the current
lines, thus reducing the quality of the reconstructed image. The two-current
approach in $2D$ is elliptic and thus does not propagate errors. The
two-current slice-by-slice reconstruction in $3D$ is also possible, but the use
of three currents seem to produce better results.

\end{enumerate}

The results of this work were presented at the conferences ``Integral Geometry
and Tomography'', Stockholm, Sweden, August 2008; ``Mathematical Methods in
Emerging Modalities of Medical Imaging'', BIRS, Banff, Canada, October, 2009;
``Inverse Transport Theory and Tomography'', BIRS, Banff, May 2010;
``Mathematics and Algorithms in Tomography'' Oberwolfach (April 2010), and
``Inverse problems and applications'', MSRI, Berkeley, August 2010. The brief
reports have appeared in \cite{KuKuAET,Oberw}.
\begin{table}[t]
\begin{center}
\begin{tabular}{|r|r|r|r|r|r|}
\hline
$j$ & $x_{j,1}$ &  $x_{j,2}$ &  $r^{\mathrm{out}}_j$ &   $r^{\mathrm{in}}_j$ &  $\alpha_j$ \\ \hline
$1 $& $-0.54$ &$ 0.54 $&$0.26 $&$0.24$ &$ 1 $  \\
$2 $& $ 0.00$ &$ 0.60 $&$0.24 $&$0.22$ &$-1 $  \\
$3 $& $ 0.60$ &$ 0.60 $&$0.16 $&$0.14$ &$ 1 $  \\
$4 $& $-0.60$ &$ 0.00 $&$0.16 $&$0.14$ &$-1 $  \\
$5 $& $ 0.60$ &$ 0.00 $&$0.26 $&$0.24$ &$-1 $  \\
$6 $& $-0.54$ &$-0.54 $&$0.26 $&$0.24$ &$ 1 $  \\
$7 $& $ 0.00$ &$-0.60 $&$0.24 $&$0.22$ &$-1 $  \\
$8 $& $ 0.60$ &$-0.60 $&$0.16 $&$0.14$ &$ 1 $  \\
$9 $& $ 0.18$ &$ 0.18 $&$0.16 $&$0.14$ &$-1 $  \\
$10$& $ 0.18$ &$-0.18 $&$0.16 $&$0.14$ &$ 1 $  \\
$11$& $-0.18$ &$ 0.18 $&$0.16 $&$0.14$ &$ 1 $  \\
$12$& $-0.18$ &$-0.18 $&$0.16 $&$0.14$ &$-1 $  \\ \hline
\end{tabular}
\end{center}
\caption{Parameters of the $2D$ phantom}
\label{T:2d}
\end{table}

\section*{Acknowledgments}

The work of both authors was partially supported by the NSF DMS grant 0908208;
the manuscript was written while they were visiting MSRI. The work of P.~K.
was also partially supported by the NSF DMS grant 0604778 and by the Award No.
KUS-C1-016-04, made to IAMCS by King Abdullah University of Science and Technology
(KAUST). The authors express their gratitude to NSF, MSRI, KAUST, and IAMCS for
the support. Thanks also go to G.~Bal, E.~Bonnetier, J. McLaughlin, L. V. Ngueyn,
L.~Wang, and Y.~Xu for helpful discussions and references. Finally, we are
grateful to the referees for suggestions and comments that helped to significantly
improve the manuscript.

\begin{table}[t]
\begin{center}
\begin{tabular}{|r|r|r|r|r|r|r|}
\hline
$j$ & $x_{j,1}$ &  $x_{j,2}$ &  $x_{j,3}$ & $r^{\mathrm{out}}_j$ &   $r^{\mathrm{in}}_j$ &  $\alpha_j$ \\ \hline
$1 $&$-0.615 $&$   -0.54  $&$   0     $&$0.26$&$   0.22 $&$   0.5  $  \\
$2 $&$-0.6   $&$   0      $&$   0     $&$0.24$&$   0.20 $&$   1    $  \\
$3 $&$0.6    $&$   0.6    $&$   0     $&$0.16$&$   0.12 $&$   0.5  $  \\
$4 $&$0      $&$   -0.6   $&$   0     $&$0.16$&$   0.12 $&$   1    $  \\
$5 $&$0      $&$   0.6    $&$   0     $&$0.26$&$   0.22 $&$   1    $  \\
$6 $&$-0.54  $&$   -0.54  $&$   0     $&$0.26$&$   0.22 $&$   0.5  $  \\
$7 $&$-0.6   $&$   0      $&$   0     $&$0.24$&$   0.20 $&$   1    $  \\
$8 $&$-0.6   $&$   0.6    $&$   0     $&$0.16$&$   0.12 $&$   0.5  $  \\
$9 $&$0.18   $&$   0.18   $&$   0     $&$0.16$&$   0.12 $&$   1    $  \\
$10$&$-0.18  $&$   0.18   $&$   0     $&$0.16$&$   0.12 $&$   0.5  $  \\
$11$&$0.18   $&$   -0.18  $&$   0     $&$0.16$&$   0.12 $&$   0.5  $  \\
$12$&$-0.18  $&$   -0.18  $&$   0     $&$0.16$&$   0.12 $&$   1    $  \\
$13$&$0      $&$   0      $&$   0.6   $&$0.18$&$   0.14 $&$   -1   $  \\
$14$&$0      $&$   0      $&$   0.6   $&$0.30$&$   0.26 $&$   1    $  \\
$15$&$0      $&$   0      $&$   -0.46 $&$0.38$&$   0.34 $&$   0.5  $  \\
$16$&$0      $&$   0      $&$   -0.46 $&$0.16$&$   0.12 $&$   0.5  $  \\   \hline
\end{tabular}
\end{center}
\caption{Parameters of the $3D$ phantom}
\label{T:3d}
\end{table}

\section*{Appendix}

In order to make it easier for the reader to repeat our simulations
we summarize in this section the details of some of our numerical experiments.

In the first two of the $2D$ simulations described in Section~\ref{S:numerical}
 we use a
 $2D$ phantom in the form of a linear combination of twelve smoothed
characteristic
functions of disks with radii $r_{j}^{\mathrm{in}}$ and centers $x_{j}$:
\begin{equation}
f(x)=\sum_{j=1}^{12}\alpha_{j}h(|x-x_{j}|,r_{j}^{\mathrm{in}},r_{j}
^{\mathrm{out}}),\quad x_{j}=(x_{j,1},x_{j,2}), \notag
\end{equation}
where
\begin{equation}
h(r,r_{j}^{\mathrm{in}},r_{j}^{\mathrm{out}})=\left\{
\begin{array}
[c]{ccc}
1 & , & r\leq r_{j}^{\mathrm{in}}\\
0 & , & r\geq r_{j}^{\mathrm{out}}\\
\exp\left[  2\frac{r_{j}^{\mathrm{out}}-r_{j}^{\mathrm{in}}}{r-r_{j}%
^{\mathrm{out}}}\exp\left(  \frac{r_{j}^{\mathrm{out}}-r_{j}^{\mathrm{in}}%
}{r_{j}^{\mathrm{in}}-r}\right)  \right]   & , & r_{j}^{\mathrm{in}}%
<r<r_{j}^{\mathrm{out}}%
\end{array}
\right.  ,     \notag
\end{equation}
and values of $\alpha_{j},$ $x_{j,1},$ $x_{j,2},$ $r_{j}^{\mathrm{in}},$ and
$r_{j}^{\mathrm{out}}$ are given in Table~1.
All the smoothed disks lie within the square computational
domain~$[-1,1] \times [-1,1]$. The forward problem was computed on
a fine $513 \times 513$ grid. We simulated propagating spherical fronts
generated by $256$ transducers equally spaced on the circle of
radius $1.6$ centered at the origin. For each transducer we simulated $257$
spherical fronts of varying radii. The reconstruction was performed on the
coarser $129 \times 129$ computational grid, from the data corresponding
to two currents. In the first experiment we used
the noiseless data, in the second one we added to the simulated values
$M_{I,J}(t,z)$ values of a random variable modeling the noise of intensity $50\%$
of the signal in $L^2$ norm. The results of these simulations are described
in Section~\ref{S:numerical}.

In Section~\ref{S:3d} we utilized a $3D$ phantom represented by a
linear combination of sixteen smoothed characteristic
functions of balls with radii $r_{j}^{\mathrm{in}}$ and centers $x_{j}$:
\begin{equation}
f(x)=\sum_{j=1}^{16}\alpha_{j}h(|x-x_{j}|,r_{j}^{\mathrm{in}},r_{j}
^{\mathrm{out}}),\quad x_{j}=(x_{j,1},x_{j,2},x_{j,3}), \notag;
\end{equation}
the values of $\alpha_{j},$ $x_{j,1}$, $x_{j,2}$, $x_{j,3}$, $r_{j}^{\mathrm{in}}$, and
$r_{j}^{\mathrm{out}}$ are given in Table~2.
In our $3D$ simulations we had to assume that the values $M_{i,j}(x)$ are
known. We modeled these values by using the above-mentioned phantom, in combination
with three boundary current profiles. In the case of the constant conductivity
these boundary currents would produce potentials equal to $x_j$, $j=1,2,3$.
We modeled the direct problem using $257 \times 257 \times 257$ computational
grid corresponding to the cube $[-1,1] \times [-1,1] \times [-1,1]$.
In the first of our $3D$ experiments the reconstruction was done on the
same grid from the noiseless data. In the second experiment the reconstruction
was done on a coarser $129 \times 129 \times 129$ grid from the data
contaminated by a $10\%$ noise (in $L^2$ norm). The results of these
reconstructions are described in Section~\ref{S:3d}.

\end{document}